\documentclass{elsarticle}

\usepackage{lineno,hyperref}
\usepackage{a4wide}
\usepackage{amsmath}
\usepackage{amsthm}
\usepackage{amssymb}
\usepackage{amsfonts}
\usepackage{mathrsfs}
\usepackage{gensymb}
\usepackage{nicefrac}
\usepackage{stmaryrd}
\usepackage{booktabs}
\usepackage{wrapfig}
\usepackage{graphicx}
\usepackage{epstopdf}
\usepackage{paralist}
\usepackage{subcaption}
\usepackage{xcolor}
\usepackage{enumitem}
\usepackage{siunitx}
\usepackage{multirow}
\usepackage{xspace}

\usepackage[capitalise]{cleveref}
\usepackage{todonotes}
\usepackage{bm}

\newcommand{\vel}{\mathbf{v}}
\newcommand{\pos}{\scalebox{0.7}{$(\mathbf{x})$}}

\newcommand{\Vavg}[1]{\left< #1 \right>} 
\newcommand{\pert}[1]{\tilde{#1}} 
 
\newcommand{\Dumux}{{Du\-Mu$^\text{x}$}}
\newcommand{\Dune}{\textsc{Dune}\xspace}

\makeatletter
\def\ps@pprintTitle{%
	\let\@oddhead\@empty
	\let\@evenhead\@empty
	\def\@oddfoot{}%
	\let\@evenfoot\@oddfoot}
\makeatother

\begin{document}
	
\begin{frontmatter}
	
	\title{Data-Driven Closure Parametrizations with Metrics: Dispersive Transport}
	
	\author[addressIWS]{Edward Coltman\corref{mycorrespondingauthor}}
	\ead{edward.coltman@iws.uni-stuttgart.de}
	\author[addressIWS]{Martin Schneider}
	\ead{martin.schneider@iws.uni-stuttgart.de}
	\author[addressIWS]{Rainer Helmig}
	\ead{rainer.helmig@iws.uni-stuttgart.de}
			
	\cortext[mycorrespondingauthor]{Corresponding author}
	
	\address[addressIWS]{Institute for Modelling Hydraulic and Environmental Systems, University of Stuttgart, Pfaffenwaldring 61, 70569 Stuttgart, Germany}
	\begin{abstract}
		This work presents a data-driven framework for multi-scale parametrization of velocity-dependent dispersive transport in porous media. Pore-scale flow and transport simulations are conducted on periodic pore geometries, and volume-averaging is used to isolate dispersive transport, producing parameters for the dispersive closure term at the Representative Elementary Volume (REV) scale. After validation on unit cells with symmetric and asymmetric geometries, a convolutional neural network (CNN) is trained to predict dispersivity directly from pore-geometry images. Descriptive metrics are also introduced to better understand the parameter space and are used to build a neural network that predicts dispersivity based solely on these metrics. While the models predict longitudinal dispersivity well, transversal dispersivity remains difficult to capture, likely requiring more advanced models to fully describe pore-scale transversal dynamics.
	\end{abstract}
\end{frontmatter}
	
\section{Introduction}	\label{sec:introduction}

	Evaluations on the Representative Elementary Volume (REV) scale are ubiquitous in the field of environmental modelling and porous medium flow research. While these REV-scale models provide efficient and generally accurate predictions, some important pore-scale phenomena that influence system dynamics are challenging to represent at these averaged scales. To solve this issue, parametrizations of multi-scale closure terms can effectively integrate pore-scale dynamics into REV-scale models, enhancing their predictive capabilities and overall reliability.	
	
	One example of a multi-scale closure term, common in the study of solute transport in porous-media flows, is dispersive transport. Transport phenomena observed at the pore scale, but not resolved by an averaged scale advection-diffusion model, are included using a dispersion term, parameterized with a tensor \cite{bear1988a}. A few models exist to define the content of this tensor, but the commonly used Scheidegger model uses the average flow velocity, along with longitudinal and tangential dispersivity parameters, $\alpha_L$ and $\alpha_T$ \cite{scheidegger1961a}. This work investigates the definition of these dispersivity parameters based on pore geometry.
	
	The use of machine learning tools is becoming increasingly common within the field of porous media research \cite{yang2024a}. In many cases, data-driven models have been introduced to improve the quality and usability of pore-scale data \cite{wang2021a}, or to help connect concepts at disparate scales \cite{tahmasebi2020a}. As is done in this work, data-driven models have also been introduced to efficiently address multi-scale closure problems \cite{lasseux2021a}. Considerable work has been done investigating the permeability parameter using this method \cite{gaerttner2021a, fu2023, mostaghimi2006a}, with further studies evaluating transport and multi-phase flow parameters \cite{graczyk2023, wang2021a, liu2019a}.  Similar efforts have begun to assess transport parameters in porous media, analyzing both diffusive transport \cite{graczyk2023} and dispersive transport \cite{kamrava2021a, taghizadeh2022a, meng2024a}. 
	This study introduces and evaluates the use of two data-driven models, an image-based model and a metrics-based model, in order to predict dispersivity parameters based on pore-geometry information.
	
	Although data-driven models are quite powerful and widely used, they are often difficult to interpret, and are not known to function well outside of the parameter space upon which they are trained \cite{kriesel2007a}. For image-based models, understanding the boundaries of this parameter space can be a challenge. To manage this, descriptions of the images used can be measured, and their distribution analyzed. In this work, descriptive metrics are developed to characterize the pore geometries. The distribution of these metrics, along with their correlation to dispersivity parameters, is then reviewed.
		
	This study introduces a data-driven framework for the multi-scale parametrization of velocity-dependent dispersive transport in porous media. For any readers primarily interested in the data-driven model development, it is recommended to move to \cref{sec:datadrivenmodel}. To begin, pore-scale flow and transport simulations are conducted on various periodic pore geometries in \cref{sec:PoreScaleEvaluations}. These results are then volume-averaged to isolate the dispersive transport, which is matched to optimal parameters at the REV scale through a dispersive closure term (\cref{sec:revscaleModel}). To better understand these results, they are validated using simple unit cells with both symmetric and asymmetric geometries. In \cref{sec:datadrivenmodel}, a convolutional neural network is trained to predict dispersivity parameters directly from images of the pore geometries, replacing the pore-scale simulations and averaging with a data-driven approach. Next, in \cref{sec:MetricsEvaluation}, descriptive metrics are introduced to provide a deeper understanding of the parameter space and to establish correlations between dispersivity and other REV-scale metrics. Finally, in \cref{sec:MNN}, these metrics are further utilized to develop a neural network that predicts dispersivity based solely on the descriptive metrics.	
	
\section{Pore-Scale Evaluations} \label{sec:PoreScaleEvaluations}
	
	Paramount to the success of any data-driven modeling effort is quality training data. In this first step, pore-scale simulations are developed and solutions are stored for use as training data. For this evaluation of dispersive transport, these pore-scale evaluations take the form of direct numerical simulations of flow and transport in periodic pore geometries. 
	
	\subsection{Flow and Transport Models}
	In order to develop pore-scale flow fields within porous structures, a porous domain, $\Omega$, is partitioned into void space, $\Omega_{\text{void}}$, and solid inclusions, $\Omega_{\text{solid}}$. Assuming stationary flow conditions within the void space, the stationary Navier-Stokes equations (\cref{eq:ns}) can be solved within the discretized void space, $\Omega_{\text{void}}$, to create a flow field:
	 
	\begin{subequations}
	  	\label{eq:ns}
	  	\begin{align}
	  		\nabla \cdot ( \varrho  \vel ) &= 0, && \\
	  		\nabla \cdot \{\varrho \vel \vel^{\mathrm{T}}
	  		- \mu (\nabla \vel + (\nabla \vel)^{\mathrm{T}})
	  		+ p \mathbf{I} \} &= q, &&
	  	\end{align}
	\end{subequations}

	\noindent with $\varrho$, the density, $\vel$, the flow velocity, $\mu$, the viscosity, $p$, the pressure, $\mathbf{I}$, the identity matrix, and $q$, a momentum source term. Flow conditions within each simulation are set to creeping, with Reynolds numbers below 1e-2. 
	
	Periodic boundary conditions are applied across all external boundaries of the system. Internal boundaries, where the void space meets solid inclusions, are set as no-slip boundaries. Flow is induced in the system with a source term, $q$, applied to the momentum balance.  

	This resolves a spatial flow field, $\vel(x)$, which can be used to evaluate a transported substance. The remaining pore-scale transport equation, evaluated on the void space, $\Omega_{\text{void}}$, is as follows: 
	 
	\begin{align}\label{eq:PS-Transport}
	  	L_{\text{Pore Scale}}(c) &:=  
	  	\frac{\partial c }{\partial t} - 
	  	\nabla \cdot (\boldsymbol{D}_{m} \nabla c) +
	  	\nabla \cdot \left(\vel c \right) = 0, \qquad \forall \,  x \in \Omega_{\text{void}}.
	\end{align}

	\noindent with $c$ as the concentration of the transported substance, $\boldsymbol{D}_{m}$ as the binary molecular diffusion coefficient. Transport conditions are described with the P\'eclet number, which is below 1 in each simulation.
	
	Flow evaluation is performed on a fully periodic pore geometry, where transport progresses from the left to the right of the system. To capture this efficiently, the domain is extended to three adjacent pore geometries in the longitudinal direction (one on each side) while maintaining periodicity tangentially. Transport progresses from the left to the right, while tangential transport continues across periodic boundaries. The inlet concentration is fixed, and outflow boundary conditions are applied on the right side. The evaluation period begins when the concentration front reaches the center pore and ends after it passes.

	\subsection{Numerical Methods}
	
	In order to solve the model described above, a finite volume Marker-and-Cell scheme can be used to discretize the domain for the evaluation of the flow field \cite{harlow1965a}. To evaluate the transport, a finite volume scheme with two-point flux approximation is used \cite{eymard2000a, versteeg2007a}. The flow and transport simulation environment \Dumux \cite{koch2020a}, based on the numerics environment \Dune \cite{bastian2008b,bastian2008a}, is used. UMFPACK \cite{davis2004a}, is used as a direct linear solver. 
	
	In addition to the standard \Dune modules, \Dune-SubGrid \cite{graeser2009a}, is used to create the computational grids for the void-space, $\Omega_{\text{void}}$, from a base structured \Dune-SPGrid grid implementation featuring periodic boundary conditions \cite{nolte2011a}. 
	
	In order to manage many pore-scale simulations and their results, simulation results and metadata are collected and stored using \texttt{JSON} based metadata collection, implemented in \cite{lohmann2022a}, and adapted for \Dumux. Some benefits of this method are the flexible type storage, and the ease of use in cross-platform processing steps. 
	
	All of the geometries, implementations, data, and models seen in this work are collected and available in a \Dumux-Pub module \cite{dumux-pub}. This module is based on the 3.7 release of \Dumux, \cite{dumux37}.
	
	\subsection{Pore Geometries and Preprocessing} \label{subsec:poregeometries}
	
	Another key aspect of pore-scale analysis is defining a pore geometry, setting the simulation domain, $\Omega_{\text{void}}$. Advances in micro-CT scanning have made digital scans of real pore geometries more accessible \cite{prodanovic2015a}, but pseudo pore geometries are also widely used due to their simplicity and design flexibility. 
	
	In this work, periodic pseudo pore geometries are created using generators based on those provided in the PoreSpy project \cite{gostick2019a}. A voronoi polygon generator, with a custom padding for periodic connectivity \cite{fritzen2008a}, and two periodic custom generators, inspired by Perlin and fractal noise tools \cite{vigier2018a}, are used within this work. Porosity, aperture sizes, scales, and filters are randomly selected in each generator to create the dataset described later in \cref{sec:datadrivenmodel}. Descriptions of the complexity of the images in this dataset are described in \cref{sec:MetricsEvaluation}.
	
	\begin{figure} 
		\centering

		\begin{subfigure}{0.31\textwidth}
			\centering
			\includegraphics[width=0.90\textwidth]{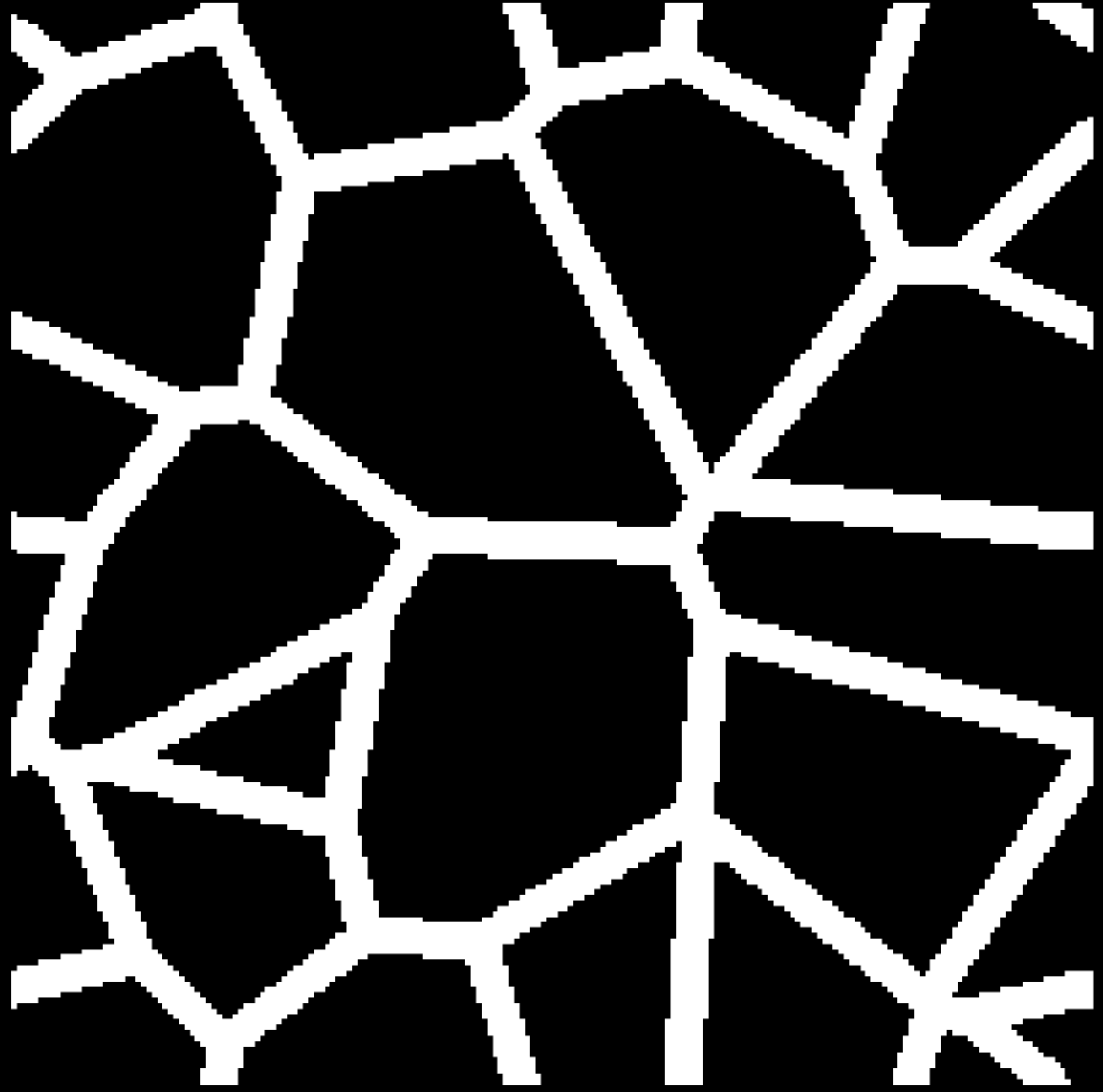}
			\caption{}
			\label{fig:unitcell_voronoi}
		\end{subfigure}
		\begin{subfigure}{0.31\textwidth}
			\centering
			\includegraphics[width=0.90\textwidth]{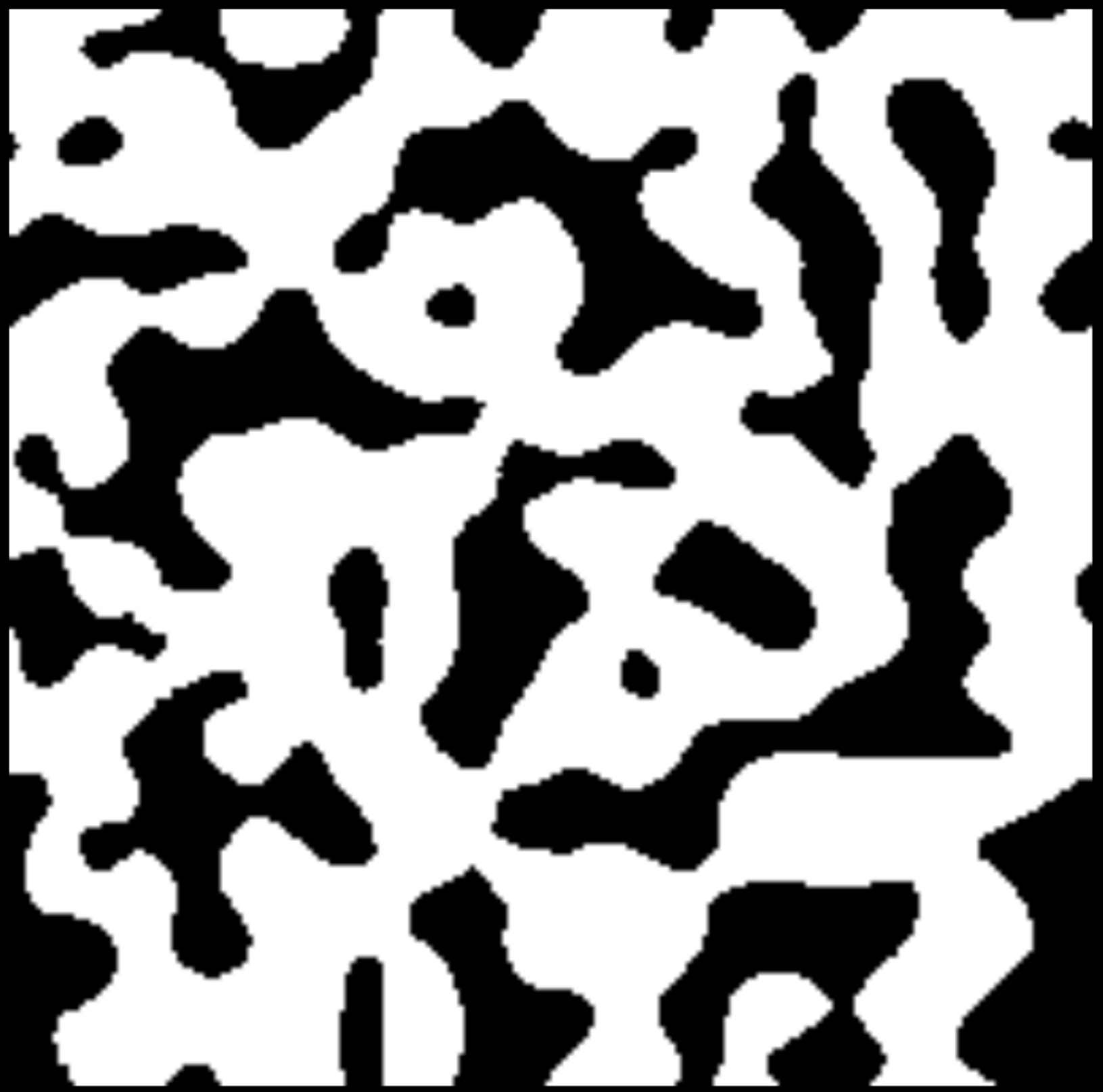}
			\caption{}
			\label{fig:unitcell_perlin}
		\end{subfigure}
		\begin{subfigure}{0.31\textwidth}
			\centering
			\includegraphics[width=0.90\textwidth]{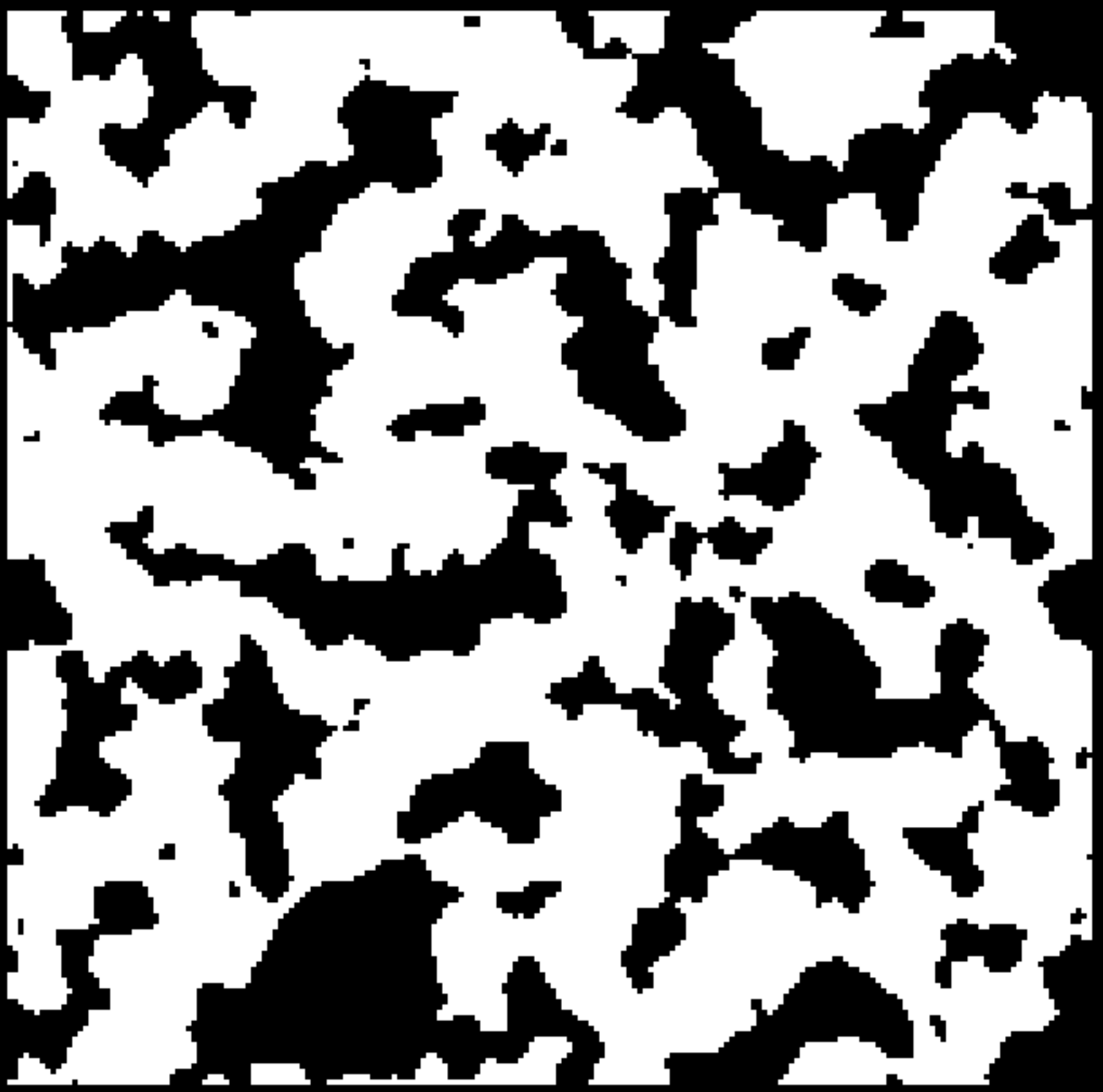}
			\caption{}
			\label{fig:unitcell_fractal}
		\end{subfigure}
		\caption{\label{fig:unitcells_complex} Examples of the pseudo pore geometries used are shown left to right: Voronoi polygons, Perlin noise, and Fractal noise based topographies. White regions belong to the void space $\Omega_{\text{void}}$, and black regions belong to the solid inclusions $\Omega_{\text{solid}}$.}
	\end{figure}
	
	To ensure usability, each unit cell is checked for connectivity, including periodicity, using a graphical component traverse \cite{berge1985a}. Disconnected pore spaces are excluded, and the image is updated. These images distinguish solid inclusions (black) from void space (white), as seen in \cref{fig:unitcells_complex}, and are stored as $200\times200$ pixel \texttt{.pbm} binary files. While 3D geometries and simulations are possible, the added complexity of a third dimension complicates interpretation. Therefore, 2D evaluations are employed in this work as a more straightforward starting point.
	
\section{REV-Scale Model Parametrization} \label{sec:revscaleModel}
	
	While pore-scale evaluations of flow and transport are valuable, they don't scale well for application-level investigations. For larger scales, evaluations at the Representative Elementary Volume (REV) scale are common in environmental modeling and porous media research \cite{bear1988a,dentz2007a,sahimi1986a}. At this scale, contaminant plume shapes and evolution are often used to parameterize models \cite{cirpka2006a, cirpka2003a}, with dispersivity calculated as a fitting parameter. As this dynamic comes from pore-scale phenomena, incorporating this into REV models requires appropriate parametrizations.
	
	\subsection{Upscaling Dispersive Transport}
	
	To transfer from pore-scale to REV-scale models, upscaling can be performed using volume averaging \cite{quintard1993a}. This method transfers high-resolution pore-scale variables to REV-scale variables, using the following decomposition:
	
	\begin{equation}\label{eq:decompose}
		\underbrace{m}_{\text{Pore-Scale Field}} = 
		\underbrace{\Vavg{m}}_{\text{Spatial Average}}
		+ \underbrace{\pert{m},}_{\text{Spatial Deviation}}
	\end{equation}
	
	where the $\Vavg{ }$ and $\pert{ }$ operators represent a spatial average and deviation, respectively. The averaged pore-scale transport equation (\cref{eq:PS-Transport}), as developed in \cite{wood2003a}, can then be given as: 
		 
	\begin{align}\label{eq:AVGPS-Transport}
		\Vavg{ L_{\text{Pore Scale}}(c)} &=  
		\frac{\partial \Vavg{c} }{\partial t} - 
		\nabla \cdot (\boldsymbol{D}_{m} \nabla \Vavg{c}) +
		\nabla \cdot \left(\Vavg{\vel} \Vavg{c}\right) 
		+ \nabla \cdot \underbrace{\left(\Vavg{\pert{\vel} \pert{c}}\right)}_{\text{Perturbation product}} = 0. &&
	\end{align}
	
	Due to the scale variance in the advective term, an additional term appears $( \nabla \cdot \left(\Vavg{\pert{\vel} \pert{c}}\right) )$. This product of the pore-scale deviations, the perturbation product, cannot be directly incorporated into an averaged REV-scale model as the deviation fields are not available at that scale. To account for its effects, a closure term, the dispersion term $(-\nabla \cdot (\boldsymbol{D}_{disp} \nabla c))$, with a dispersion tensor $\boldsymbol{D}_{disp}$ (simplified here as $\boldsymbol{D}$), is added to capture directional transport. For further details on the assumptions made in the development of \cref{eq:AVGPS-Transport}, please refer to \cite{whitaker1999a, wood2003a}. 
		
	\subsection{Dispersive Transport Model}
	
	At the REV scale, the transport equation takes a similar form to that of the averaged pore-scale equation, with a dispersion term replacing the perturbation product term. This can be written as follows:	
	
	\begin{align} \label{eq:transportREV}
		L_{\text{REV}}(\Vavg{c}) &:=  
		\frac{\partial \Vavg{c}}{\partial t} - 
		\nabla \cdot (\boldsymbol{D}_{m} \nabla \Vavg{c}) +
		\nabla \cdot ( \Vavg{\vel} \, \Vavg{c} ) -
		\nabla \cdot \underbrace{\left(\mathbf{D}(\Vavg{\vel},\mathbf{\alpha}) \nabla \Vavg{c} \right)}_{\text{Dispersion}} = 0, &&
	\end{align}

	\noindent where the diffusive and dispersive term are often combined to an effective dispersion term. The dispersion term added here requires a dispersion tensor, $\mathbf{D}$, which is defined in this work using the following relationship:
	
	\begin{equation}\label{eq:dispersionmodel}
		\mathbf{D}_{disp} \text{ or } \boldsymbol{D}(\Vavg{\vel}, \boldsymbol{\alpha}) := \qquad
		\left(\alpha_L - \alpha_T\right) \frac{\Vavg{\vel} \Vavg{\vel}^T}{\|\Vavg{\vel}\|} + 
		\alpha_T  \|\Vavg{\vel} \|\boldsymbol{I}.
	\end{equation}

	Here, the averaged velocity, $\Vavg{\vel}$, the longitudinal dispersivity, $\alpha_L$, and the transversal dispersivity, $\alpha_T$, determine the form of the tensor $\mathbf{D}$. This form of the dispersion tensor ensures that the averaged velocity vector $\Vavg{\vel}$ and its orthogonal $\Vavg{\vel}^\perp$ are eigenvectors of the dispersion tensor:
	\begin{equation}
		\mathbf{D} \Vavg{\vel} = \alpha_L \|\Vavg{\vel}\| \Vavg{\vel}\, \qquad
		\mathbf{D} \Vavg{\vel}^\perp = \alpha_T \|\Vavg{\vel}\| \Vavg{\vel}^\perp.
	\end{equation}
	
	From a numerics standpoint, this tensor can easily be discretized in a finite volume REV-scale model, where $\alpha_L$ and $\alpha_T$ are provided as parameters, and the average velocity, $\Vavg{\vel}$, reconstructed from momentum fluxes \cite{srinivasan2013a}.
	
	\subsection{Optimal Parameters}\label{subsec:OptimalAlphas}
		
	Instead of solving a closure problem \cite{wood2003a}, the pore-scale solution can be used to find the dispersivity parameters, $\alpha_L$ and $\alpha_T$. This is done by equating the pore-scale averaged perturbation product with the REV-scale dispersion term, as seen in the following minimization target:
	
	\begin{equation}\label{eq:optimaldispersion}
		\boldsymbol{\alpha} = \arg \text{min} \mathcal{L}_2 (\boldsymbol{\alpha}) , \qquad 
		\mathcal{L}_2(\boldsymbol{\alpha}) := 
		\left| \left| 
		\Vavg{\pert{\vel} \pert{c}} - 
		\boldsymbol{D}(\Vavg{\vel}, \boldsymbol{\alpha}) \left(\text{-} \nabla \Vavg{c}\right) 
		\right| \right|^2
	\end{equation}
	
	From the pore-scale analysis, all vectors observed in this minimization target, the averaged velocity $\Vavg{\vel}$, the negative concentration gradient, $\left(\text{-} \nabla \Vavg{c}\right)$, and the perturbation product, $\Vavg{\pert{\vel} \pert{c}}$, can be calculated from a pore-scale simulation. The perturbation product and the concentration gradient change throughout the transient pore-scale evaluation, and are collected at multiple time intervals during each pore-scale simulation. 
	
	Using these collected data points, optimal values for the length parameters  $\alpha_L$ and $\alpha_T$ can be found for a given pore geometry. Inserting \cref{eq:dispersionmodel} yields the following equations for optimal length parameters:
	
	\begin{equation}\label{eq:alphaEquation}
		\alpha_L = 
		\frac{\left( \Vavg{\pert{\vel} \pert{c}} \cdot 
			\left(\frac{\Vavg{\vel}}{\|\Vavg{\vel}\|}\right)\right)}
		{\|\Vavg{\vel}\| \left(\text{-} \nabla \Vavg{c} 
			\left(\frac{\Vavg{\vel}}{\|\Vavg{\vel}\|}\right)\right)}
		\qquad 
		\alpha_T = 
		\frac{\left( \Vavg{\pert{\vel} \pert{c}} \cdot 
			\left(\frac{\Vavg{\vel}^\perp}{\|\Vavg{\vel}^\perp\|}\right)\right)}
		{\|\Vavg{\vel}\| \left(\text{-} \nabla \Vavg{c} 
			\left(\frac{\Vavg{\vel}^\perp}{\|\Vavg{\vel}^\perp\|}\right)\right)}.
	\end{equation}
	
	To ensure a semi-definite dispersion tensor, two constraints must be applied. First, the transversal dispersivity must be non-negative, $\alpha_T \geq 0$. Second, the transversal dispersivity should not exceed the longitudinal dispersivity, represented as $\alpha_T \leq \alpha_L$. With these constraints in place, parametrizations can be determined for each evaluated time step. An average of the dispersivity parameters calculated for each time step can then be taken as the optimal dispersion parametrization for the specified pore geometry. 

	The dispersion model used here, defined using the average velocity vector and two dispersivity parameters, is a simplification of the full dispersion dynamics observed at the pore scale. This model can only create a fully symmetrical tensor and assumes a primary longitudinal axis aligned with the velocity \cite{bear1988a}. Limitations of this approach are discussed in \cite{lichtner2002a}, which investigates an extended model incorporating four characteristic lengths ($\alpha_{1-4}$), the average velocity vector, $\Vavg{\vel}$, and a geometric symmetry axis vector, $\boldsymbol{\lambda}$. Here, it is explained that the dispersion model in  \cref{eq:dispersionmodel} is only valid when the axis $\boldsymbol{\lambda}$ is either zero or aligned with $\Vavg{\vel}$; otherwise, additional $\alpha$ parameters are needed for a representative dispersion tensor. 

	In cases where the $\Vavg{\pert{\vel} \pert{c}}$ vectors do not align within the same quadrant as $\text{-} \nabla \Vavg{c}$, negative $\alpha_T$ values indicate inadequacy in the simplified model for describing transversal dispersion, which are then constrained to 0. Although the dispersive perturbation product can be isolated at the pore scale, a comprehensive dispersion model for the REV scale that clearly defines parameter dependencies is lacking \cite{lichtner2002a}. While it is possible to develop further parametrizations using volume-averaged data and more complex models, such approaches, incorporating multiple dispersivity parameters and a symmetry axis, are not commonly employed at the REV scale and would significantly increase model complexity.
	
	\subsection{Volume Averaging}
	
	In order to transfer the high-resolution pore-scale results to the REV scale, a spatial averaging is required \cite{quintard1993a}. The vectors required for the calculation of the dispersive parameters in \cref{eq:optimaldispersion,eq:alphaEquation}, can be collected by averaging the pore-scale simulation solutions. Two averaging methods are used to manage this, a convolutional averaging and a unit averaging method. 
	
	\begin{figure}
		\centering
		\begin{subfigure}{0.45\textwidth}
			\centering
			\includegraphics[width=0.95\textwidth]{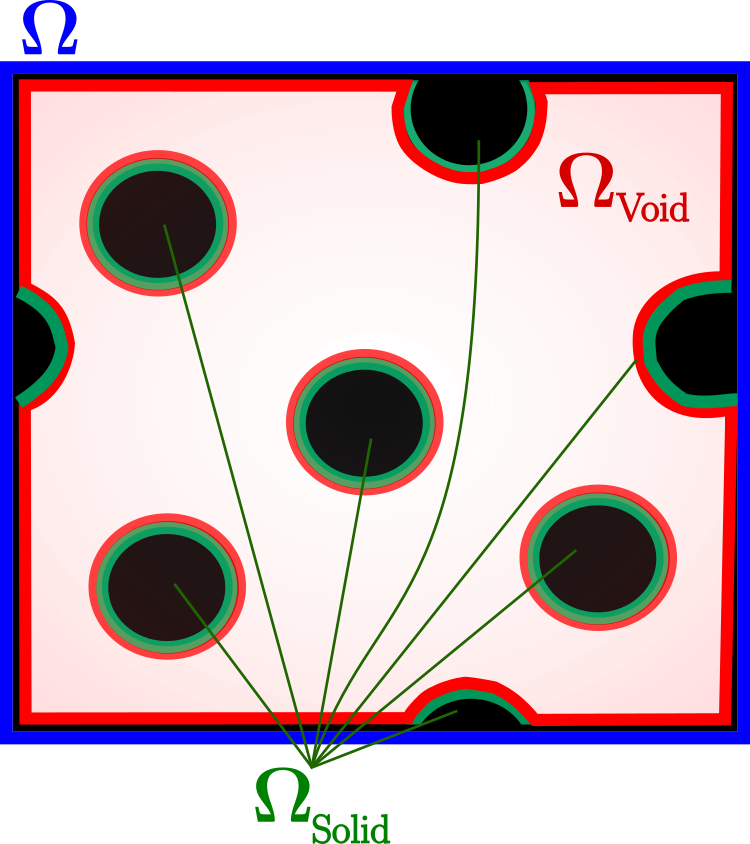} 
			\caption{\label{fig:UnitAvg}}
		\end{subfigure}
		\begin{subfigure}{0.45\textwidth}
			\centering
			\includegraphics[width=0.95\textwidth]{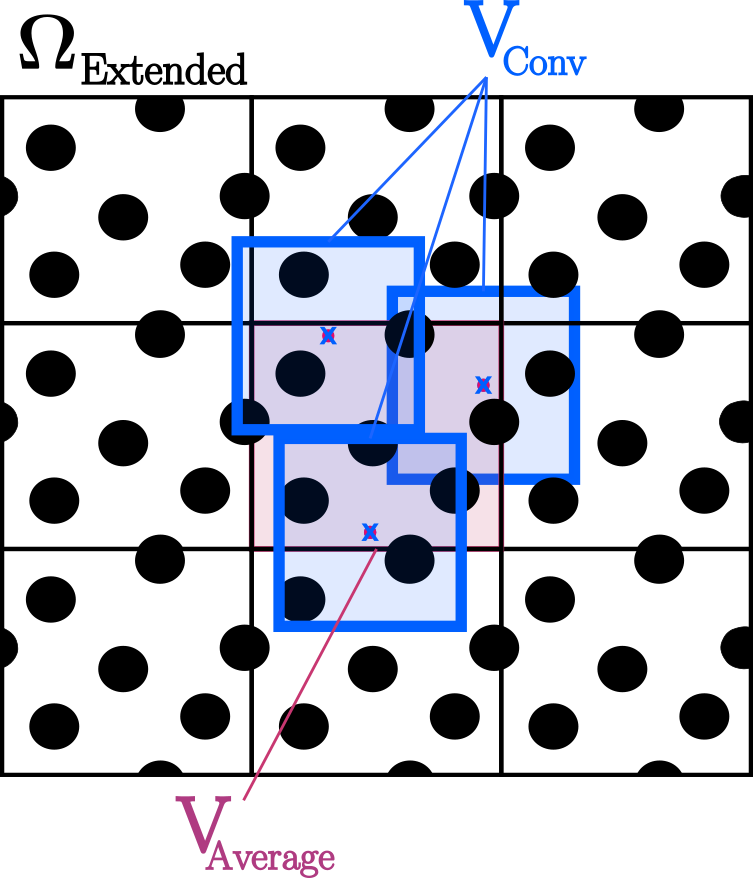}
			\caption{\label{fig:ConvAvg}}
		\end{subfigure}
		\caption{\label{fig:averaging} A depiction of the various volumes used while performing a volume averaging of flow and transport at the pore scale. In \cref{fig:UnitAvg}, the full domain, $\Omega$, and its split void-solid regions, $\Omega_{\text{void}}$ and $\Omega_{\text{solid}}$, are shown. In \cref{fig:ConvAvg}, the extended domain, $\Omega_{\text{Extended}}$, and the volumes related to convolutional and unit averaging, $\text{V}_{\text{Conv}}$ and $\text{V}_{ \text{Average}}$, are shown.} 
	\end{figure}
	
	A convolutional or ensemble average \cite{dagan1990a,dagan1991a,wood2003a}, is first used to develop an averaged flow and transport concept at the original pore-scale grid resolution. As shown in \cref{fig:ConvAvg}, flow and transport results are extended to a $3\times3$ grid of the original periodic grid geometry, $\Omega_{\text{Extended}}$. The central periodic geometry then defines the averaging volume, $\text{V}_{\text{Average}}$, and for each solution point, $\mathbf{x}$, within the averaging volume, $\text{V}_{Average}$, the following averaging is applied:
		
	\begin{equation}\label{eq:convaverage}
		\text{for all } \pos \text{ in } \text{V}_{Average} \text{:} \,\,
		\Vavg{c}\pos = 
		\frac{1}{\text{V} \pos } \int_{\text{V}} c \,d \text{V} \pos
		\,\, \text{, with } \,\, \text{V} \pos = \text{V}_{Conv}.
	\end{equation}
	
	In this work, the size of the convolutional filter $\text{V}_{\text{Conv}}$ has the same shape and size as the original domain $\Omega$, and no kernel weighting function is applied during averaging. 
	
	A unit average is then applied to condense a high-resolution field, $\mathcal{F}$ to a single representative value. This is done using the same averaging volume, $\text{V}_{\text{Average}}$, and the following equation:
		
	\begin{equation}\label{eq:unitaverage}
		\Vavg{\mathcal{F}} = \frac{1}{\text{V}}\int_{\text{V}} \mathcal{F} \,d\text{V},
		\,\, \text{ with } \,\, \text{V} = \text{V}_{Average}.
	\end{equation}
	
	In order to produce the averaged velocity vector, $\Vavg{\vel}$, a simple unit averaging can be performed on the original pore-scale flow solution. As this pore-scale field is stationary and fully periodic, the convolutional averaging step can be removed as it produces the same result. 
	
	To define the perturbation product, $\Vavg{\pert{\vel} \pert{c}}$, pore-scale deviations from the average are calculated as in \cref{eq:decompose}. First, a convolutional averaging generates a high-resolution averaged concentration field, $\Vavg{c}_{\text{Conv}}$. The deviation between the original concentration, $c$, and this average, $\Vavg{c}_{\text{Conv}}$, is isolated and multiplied by the deviation between pore-scale velocity, $\vel$, and the averaged velocity, $\Vavg{\vel}$. This yields one high-resolution perturbation product field per time step, $\pert{\vel} \pert{c}$, which is then unit averaged to yield $\Vavg{\pert{\vel} \pert{c}}$.

	To find the negative concentration gradient, $\left(\text{-} \nabla \Vavg{c}\right)$, the gradient of the high-resolution field, $\Vavg{c}_{\text{Conv}}$, is measured at each pixel and averaged to produce a concentration gradient vector for each time step.
	
	\subsection{Verification}	
	
	To evaluate the methods used to parametrize a REV-scale dispersion model from pore-scale data, a series of verification tests are conducted. Simple unit cells with solid inclusions, such as circles, squares, ellipses, and triangles, are used to observe parametrizations and compare them against prescribed shape variations. Although not representative of real porous media, these shapes are easy to describe, and their effects on flow and transport are simpler to assess. The shapes are varied in size, aspect ratio, and rotation, and as they do not intersect boundaries or contain disconnected features, no preprocessing is required. Such geometries are commonly used in pore-scale analysis \cite{strohbeck2023a} and experiments \cite{weishaupt2020a}.
	
	\subsubsection{Symmetric Cases: Circles and Squares}
			
	\begin{figure}
		\centering
		\includegraphics[width=0.96\textwidth]{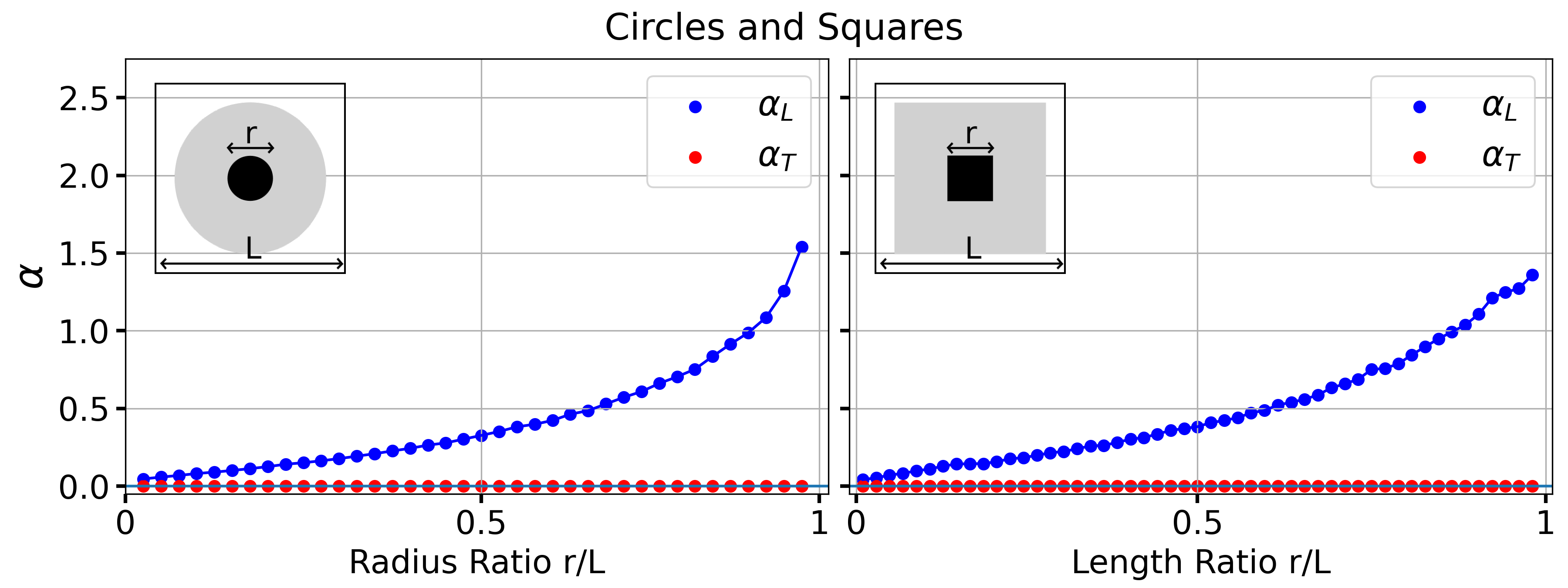}
		\caption{ The development of $\alpha_L$ and $\alpha_T$ for symmetrical cases. Considered are circular and rectangular solid inclusions of increasing length ratio ${\text{r}}/{\text{L}}$.}
		\label{fig:alphas_CircleSquare}
	\end{figure}

	The first set of verification cases deals with symmetrical solid inclusions, with either a square or circle shape centered in the domain. The radius or side length of these inclusions is increased until flow through the porous cell is nearly blocked. This progression can be observed from left to right in \cref{fig:alphas_CircleSquare}.

	As all vectors in these cases are aligned with the x-axis, due to the geometrical and flow symmetry, the transverse dispersivity, $\alpha_T$, should be zero. For the smallest solid inclusions, the longitudinal dispersivity, $\alpha_L$, approaches zero. With larger solid inclusions, the flow field disruption increases, with the highest longitudinal dispersivity values occurring in cases with the tightest bottlenecks.
	
	\subsubsection{Asymmetric Cases: Ellipses and Triangles}

	In addition to the symmetric cases, asymmetric shapes can be evaluated to compare developments in the transversal direction. The asymmetric shapes evaluated here are variations on an ellipsoidal and a triangular solid inclusion. In these cases, the averaged velocity and perturbation product vary in the transversal direction such that $\alpha_T$ may become non-zero. 
	
	\begin{figure}
		\centering
		\includegraphics[width=0.96\textwidth]{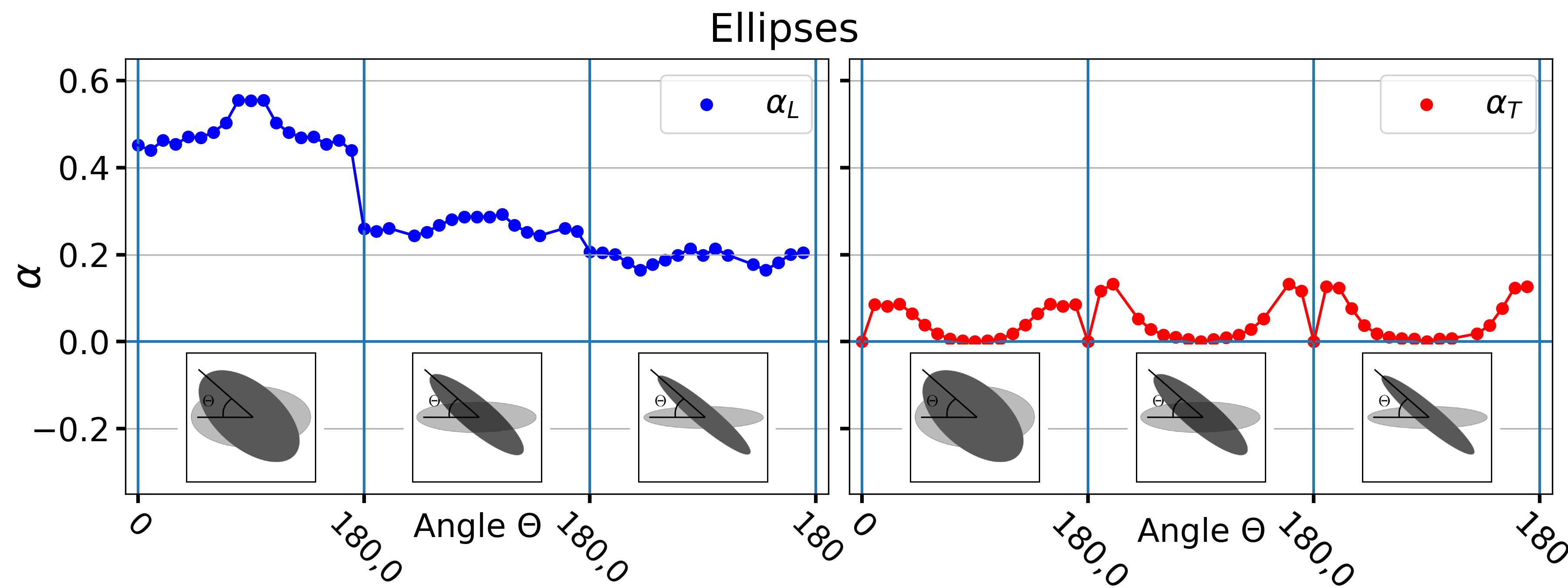}
		\caption{The development of $\alpha_L$ and $\alpha_T$ for asymmetrical ellipsoid cases. The aspect ratio and the pitch rotation $\Theta$ are varied. }
		\label{fig:alphas_ellipses}
	\end{figure}	
		
	As seen in \cref{fig:alphas_ellipses}, ellipsoids are rotated around a central position, varying in aspect ratio and rotational pitch. Within each section of  \cref{fig:alphas_ellipses}, the aspect ratio is constant, and increases from left to right. In each section the rotation pitch changes, with a symmetrical $90^\circ$ case at the center. As seen in the symmetrical cases, ellipsoids with higher aspect ratios (thinner ellipsoids) generally have smaller longitudinal dispersion, while tighter longitudinal bottlenecks (e.g., at $90^\circ$) increase $\alpha_L$. For transversal dispersion, $\alpha_T$ remains zero in symmetrical cases ($0^\circ, 90^\circ, 180^\circ$) but increases when the inclusion is slightly misaligned from the x-axis.
	
\begin{figure}
		\centering
		\includegraphics[width=0.96\textwidth]{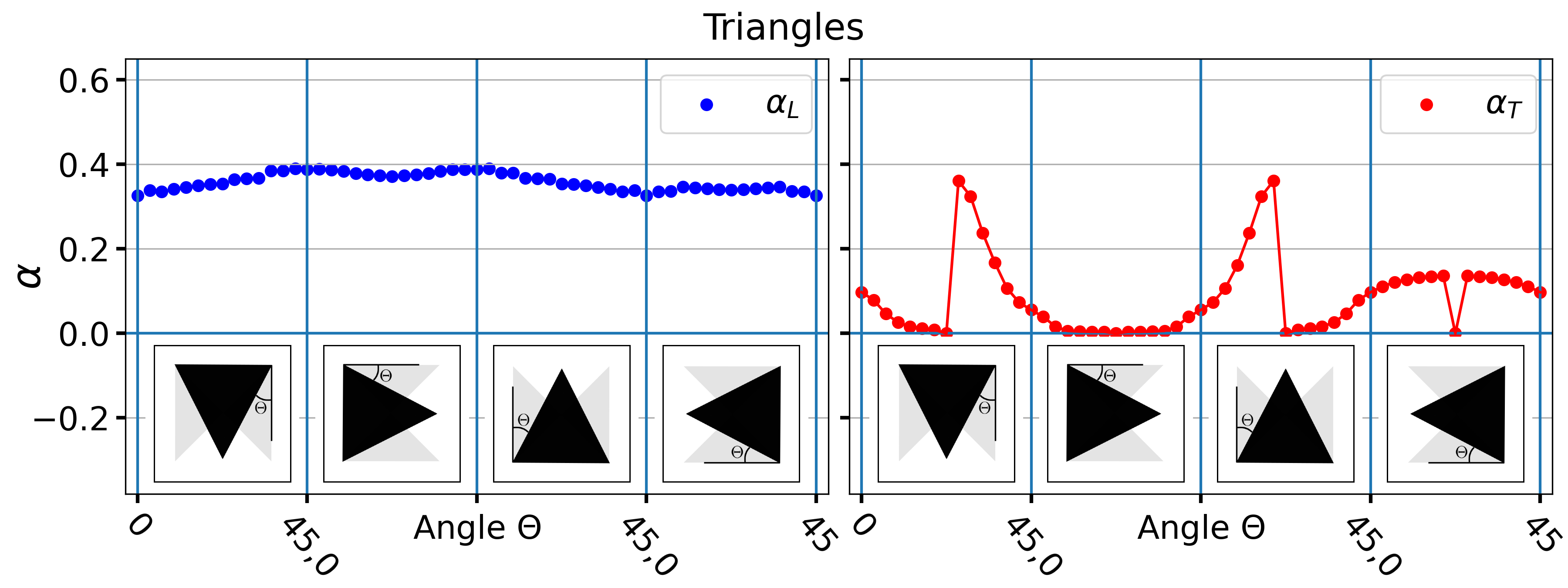}
		\caption{The development of $\alpha_L$ and $\alpha_T$ for asymmetrical triangular cases. The face orientation, $\Theta$, is varied. }
		\label{fig:alphas_triangles}
	\end{figure}

	The triangular cases introduce further asymmetry. Each triangular inclusion has a similar volume and position but varies in orientation. In \cref{fig:alphas_triangles}, these cases are divided into four sections, with one triangle face parallel to a unit cell face per section, and the remaining faces rotated between $0^\circ$ and $45^\circ$.
	
	The first and third sections, which are x-axis reflections of each other, show symmetrical trends in $\alpha_L$ and $\alpha_T$. When the triangle’s upstream face is closest to the unit cell border, $\alpha_L$ is higher. In cases of y-axis symmetry, $\alpha_T = 0$, and as with ellipses, slight deviations from symmetry result in the highest $\alpha_T$.
		
	The second and fourth sections involve y-axis reflections, and show similar trends in $\alpha_L$ due to comparable flow bottlenecks. Slight differences arise due to the varied position of the bottleneck. $\alpha_T$ is higher in the fourth section, where angled faces cause more outward transport than in the second section.	

\section{Multi-Scale Data-Driven Model}	\label{sec:datadrivenmodel}
	
	As REV-scale parameters depend on pore-scale properties, a data-driven model can be developed to replace these simulation, averaging, and scale-bridging steps outlined above. This model should directly predict REV-scale parameters from pore-scale data, represented here as pore geometry images, using a Convolutional Neural Network (CNN). CNNs, designed for spatial data like images, process inputs through layers of convolutional filters, pooling, and computational nodes with adjustable weights, $\xi$, to make predictions \cite{kriesel2007a}. These models are effective at capturing complex, multivariate relationships \cite{goodfellow2016a} and are well-suited for solving the non-linear regression problem at hand.

\subsection{Convolutional Neural Networks (CNNs)}\label{subsec:cnns}

\begin{figure}
	\centering
	\includegraphics[width=0.85\textwidth]{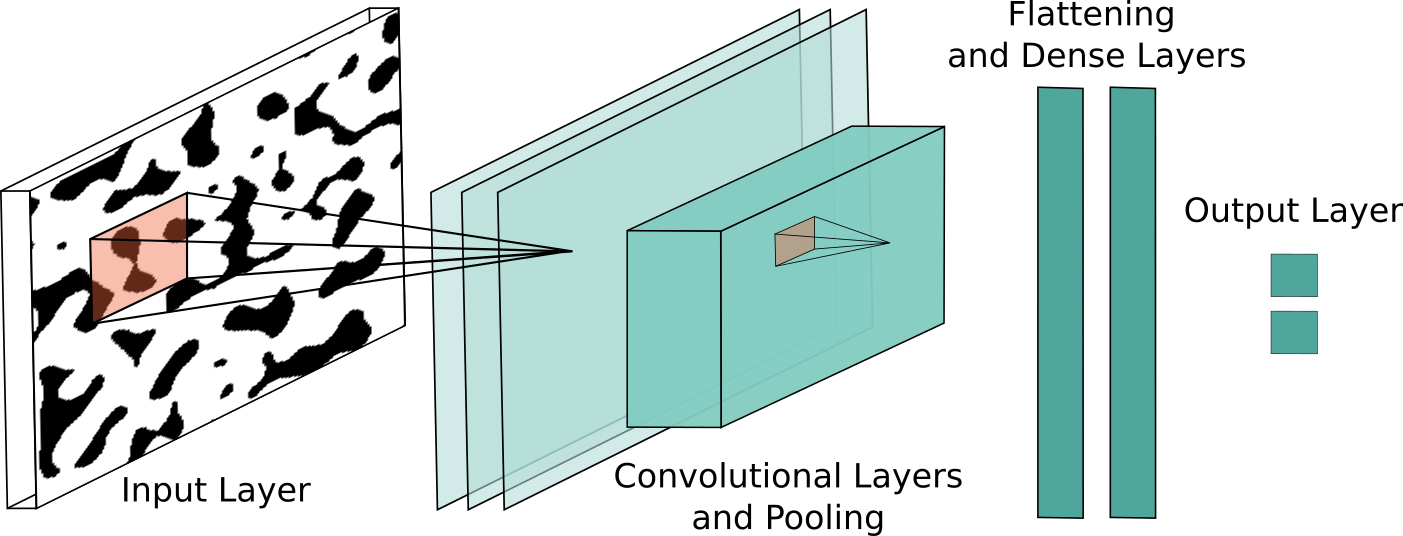}
	\caption{\label{fig:CNNs} An example visualization of a CNN (Convolutional Neural Network). Images, provided as input, are passed through a series of learned convolutional filters and pooling layers. The extractions of these layers are then flattened and passed through a series of learned dense layers, where regression output values are predicted.} 
\end{figure}

	In many cases, the parameters used at the REV scale are primarily influenced by pore-scale geometry. For this reason, pore geometries is often used as input data for data-driven models predicting these parameters \cite{gaerttner2021a, kamrava2021a, wu2019}, with various architectures tested \cite{graczyk2023}. While standard dense layers can process unordered input data, convolutional neural networks (CNNs) are more effective for classification and regression tasks with image input \cite{simonyan2014a}. By applying convolutional filtering and pooling layers, followed by dense layers leading to a regression output, a CNN can be trained to predict REV-scale parameters from pore geometry images. An example CNN is shown in \cref{fig:CNNs}.

\subsubsection{Dataset Review}

\begin{figure}
	\centering
	\includegraphics[width=0.96\textwidth]{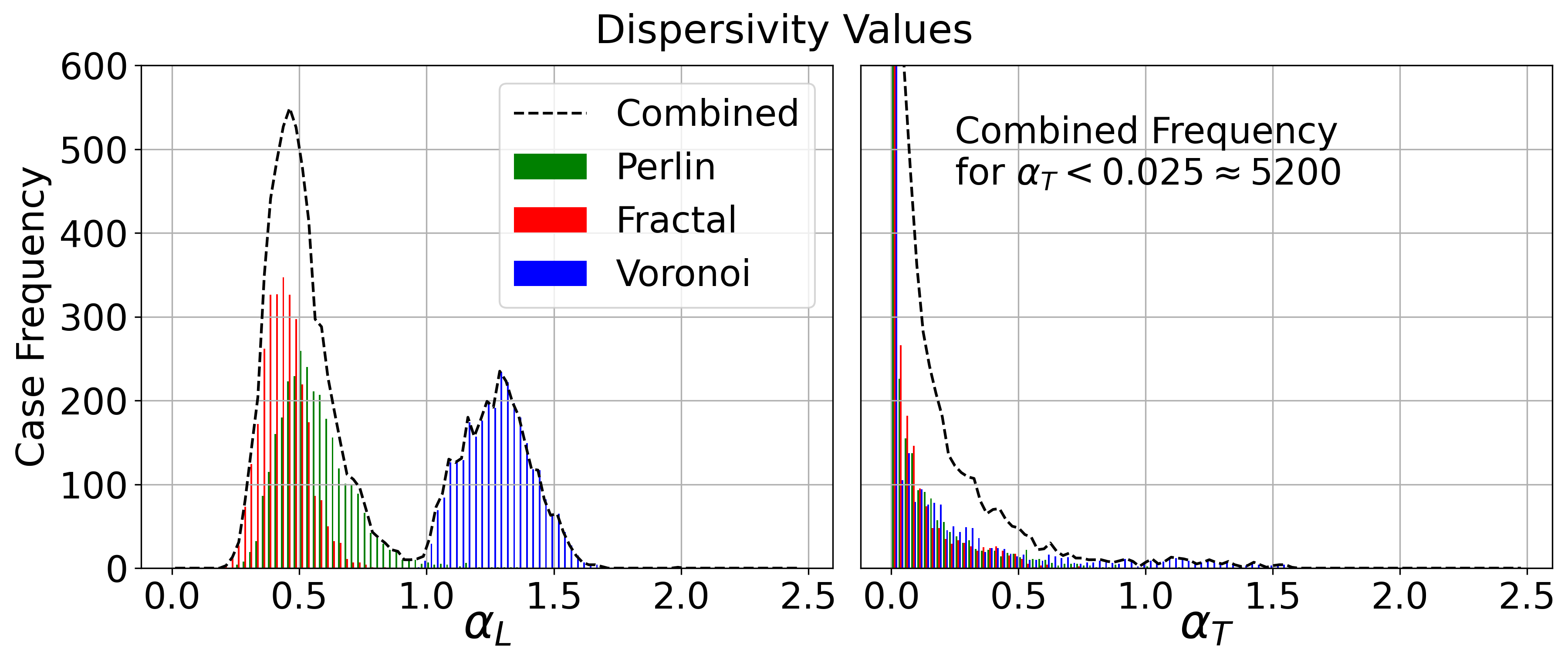}
	\caption{\label{fig:Alpha_dataset} The distribution of the dispersivities $\alpha_L$ and $\alpha_T$ developed in the dataset. }
\end{figure}

To develop a dataset upon which the data-driven model can be trained, a series of pore geometries are developed. For each type of complex pore geometry (see \cref{subsec:poregeometries} and \cref{fig:unitcells_complex}), 3000 instances are created, for a total of 9000 cases. For each case, the flow and transport problems are solved and volume averaged. With these results, the optimal dispersivity coefficients, $\alpha_L$ and $\alpha_T$, are found for each case. The distribution of the coefficients making up the dataset is shown in \cref{fig:Alpha_dataset}. 

The voronoi cells produce longitudinal dispersivities typically between 1.0 and 1.5, whereas the fractal cells focus in a lower range, from 0.2 to 0.75, and the perlin cells between 0.25 and 1.0. For the transversal dispersivity, $\alpha_T$, nearly all geometries produce values below 0.5, clustered primarily closer to 0.0. This clustering close to zero can be explained by the constraints outlined in \cref{subsec:OptimalAlphas}. For this parameter, a clear tendency for a certain type of cell is not observed. 

\subsubsection{Model Architecture}

The convolutional neural network used in this work is based on the VGG-16 architecture \cite{simonyan2014a}, with modifications to the output and input layers for model regression and to accommodate periodic 200x200 binary images. As detailed in \cref{tab:vgg16}, the architecture is divided into a feature extraction section (layers 1-18) and a dense calculation section (layers 18-21). The feature extraction includes convolutional blocks interspersed with 2x2 max pooling layers. After the last pooling layer, the results are flattened and passed to dense layers, starting with a fully connected layer of 4096 nodes. The final dense layer is split into two heads, each with 2048 nodes, allowing separate optimization for the two regression targets. All layers use the rectified linear unit (ReLU) activation function, except for the output layer, which uses a linear function for regression. In total, this model has 134 million learnable parameters ($N_\xi = 134$M) in a feed-forward network. 

One variation from standard architectures is the inclusion of a periodic specific convolutional layer at the beginning. These custom layers are introduced to the model to account for the periodicity of the pore geometry images, such that boundaries are accounted for in the same manner they are accounted for in the pore-scale simulations. This layer, based on the work presented in \cite{schubert2019a}, adds a periodic padding to each convolutional layer allowing wrap around filters in each direction, and is developed as an extraction of the standard 2D TensorFlow/Keras convolutional layer.  

\begin{table}
	\centering
	\begin{tabular}{|c|c|c|c|c|c|} \hline
		\textbf{Layer} & \textbf{Type} & \textbf{Filters $\&$ Size} & \textbf{Output Size} & \textbf{Weights $\xi$ per Layer} \\ \hline
		1-2 & 2x Convolutional Periodic & 64, $3\times3$ & $204\times204$ & 37,568\\ 
		3 & Max Pooling & $2\times2$ & $102\times102$ & 0 \\ \hline
		4-5 & 2x Convolutional & 128, $3\times3$ & $102\times102$ & 221,440 \\ 
		6 & Max Pooling &  $2\times2$ & $56\times56$ & 0 \\  \hline
		7-9 & 3x Convolutional & 256, $3\times3$ & $56\times56$ & 1,475,328 \\ 
		10 & Max Pooling & $2\times2$ & $26\times26$ & 0 \\ \hline
		11-13 & 3x Convolutional & 512, $3\times3$ & $26\times26$ & 5,899,776 \\ 
		14 & Max Pooling & $2\times2$ & $13\times13$ & 0 \\	\hline
		15-17 & 3x Convolutional & 512, $3\times3$ & $13\times13$ & 7,079,424 \\ 
		18 & Max Pooling & $2\times2$ & $7\times7$ & 0 \\ \hline
		19 & Fully Connected &  - & 4096 & 102,764,544 \\ 
		$20_{a,b}$ & 2x Split Head Dense Layers &  - & 2048 & 16,781,312 \\ \hline
		$21_{a,b}$ & 2x Regression Output &  - & 1 & 4098 \\ \hline
	\end{tabular}
	\caption{\label{tab:vgg16} Dispersivity CNN Architecture}
\end{table}

\subsubsection{Loss Function and Optimizer}

During training, the Mean Square Error (MSE) loss function with $\mathcal{L}^2$ regularization is used to find optimal model weights. This is defined here as the following:

\begin{equation} \label{eq:MSE}
	\left(\text{MSE} + \mathcal{L}^2\right) = \frac{1}{n} \sum_{i=1}^{n}(\boldsymbol{\alpha}_i - \hat{\boldsymbol{\alpha}}_i)^2 + \lambda_R \sum_{j=1}^{N_\xi} \xi_j^2, 
\end{equation}

\noindent where $n$ is the total number of training samples, and $\boldsymbol{\alpha}_i$ and $\hat{\boldsymbol{\alpha}}_i$ represent the target and predicted dispersivities. The second term represents a $\mathcal{L}^2$ regularization term, where $p$ is the total number of learned model parameter weights,  $\xi_j$, that are each regularized by the hyper parameter $\lambda_R$. While $\lambda_R$ defaults to 0, it can be increased in some layers to promote smaller weights and reduce overfitting.

The optimization algorithm used during training is the adaptive moment estimation method, known as the \texttt{'adam'} algorithm \cite{kingma2014a}. This method adjusts the learning rate based on the history of gradients and is widely adopted in the field. A learning rate of $\epsilon_l = 0.001$ is selected for this model.

\subsubsection{Dataset Splitting, Normalization, and Augmentation}

The full dataset of 9000 images and dispersivities is shuffled and divided into three subsets: $70\%$ for training, $20\%$ for validation, and $10\%$ for testing. The training set is used directly for model training, while the validation set helps determine when to stop training based on early stopping criteria. If the validation loss does not decrease after a predefined number of steps, training stops. The independent test set is reserved for final model evaluation.

To address the varied distributions of dispersivity data, the output is normalized to a range of 0 to 1, allowing for more balanced treatment in the loss function. Initially, a Box-Cox transformation \cite{boxcox1964a} is applied to make the data more normally distributed, followed by linear scaling to fit the 0 to 1 range. The normalized predictions can later be reconstructed to their original form using the inverse transformations. With this, the multi-regression predictions are more equally treated in the loss function, to ensure balanced training across both variables.
	
\begin{figure}
	\centering
	\includegraphics[width=0.70\textwidth]{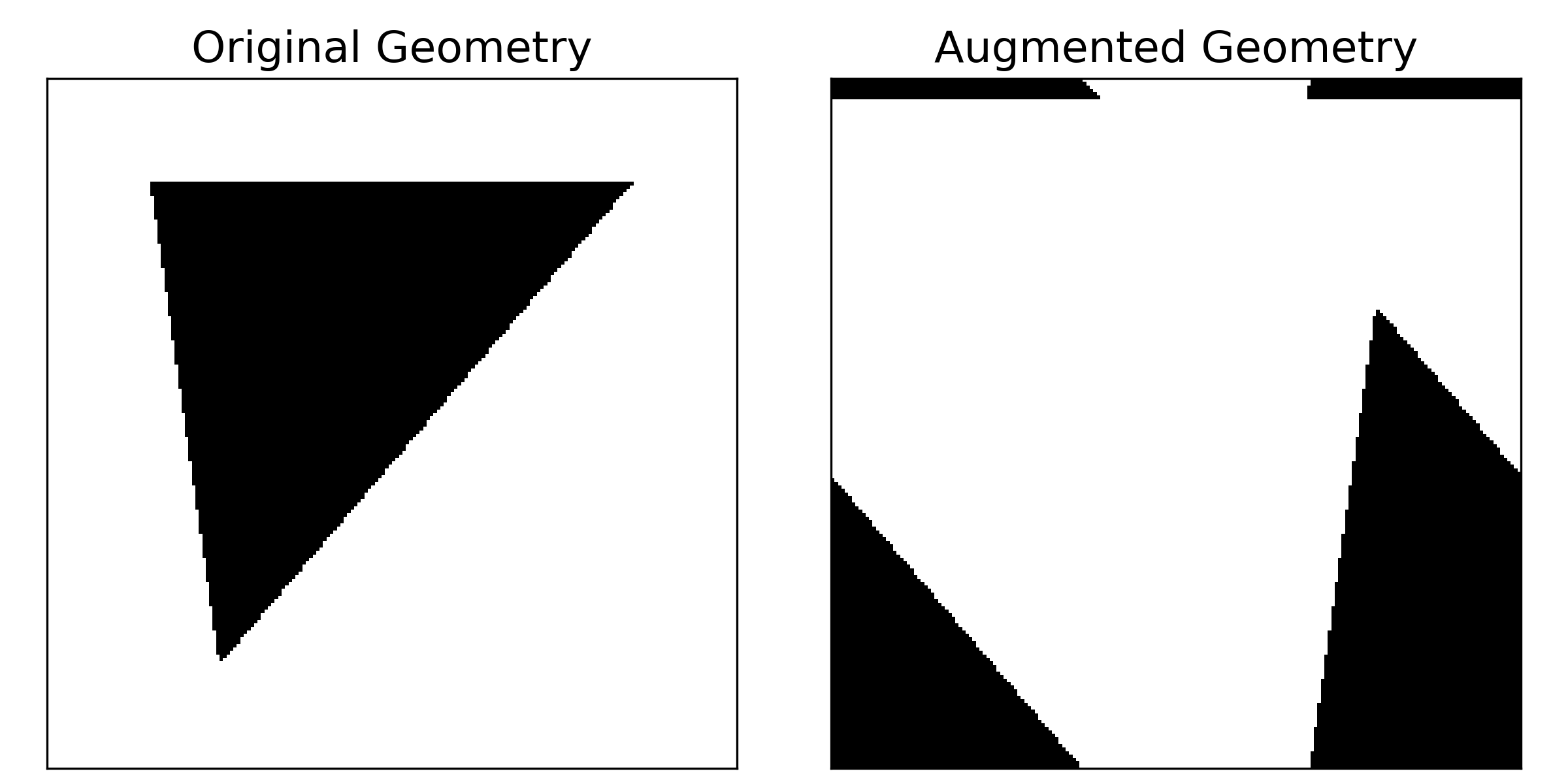}
	\caption{\label{fig:augmentation} Shown is an example of an augmented unit cell, both flipped around the x-axis, and shifted periodically. }
\end{figure}

To maximize the training set's content, data augmentation techniques are employed, including flipping images around the x-axis and applying periodic shifts. This approach can double the training set size with flipped images and increase it further through random shifts. Data augmentation occurs randomly during each training epoch, ensuring batches are altered in each iteration, which helps reduce the risk of overfitting while managing storage and memory constraints. An example of this augmentation is shown in \cref{fig:augmentation}. Here, one of the simple triangle cases is shown to more clearly see the changes made, with a flipping around the x-axis and a periodic vertical and horizontal shift. 

\subsubsection{Implementation and Availability}

All neural networks built within this framework use the toolkit TensorFlow \cite{tensorflow2015a}. Layers, connectivity, loss functions and data management tools are all available within this toolkit, and models can be built as a modular collection of these building blocks. As previously mentioned, customizations are developed for a periodic convolutional layer, a dataset generator for augmented images, as well as normalization and reconstruction layers for the output data. 

Training these models can be performed on a local machine using a training script provided in \cite{dumux-pub}. Additionally, training can be performed using the provided jupyter notebook and dataset using the Google Colaboratory cloud computing platform \cite{bisong2019a}. With this, the reproducibility of the training is easily possible, and high-end run times can be developed easily on the cloud based GPU hardware available. These scripts and notebooks, as well as saved model weights and architectures, are all available within the \texttt{machinelearning} directory in the \texttt{dumux-pub} repository mentioned above \cite{dumux-pub}. 

\subsection{CNN Model Performance}\label{subsec:CNNPerformance}
\begin{figure}
	\centering
	\includegraphics[width=0.96\textwidth]{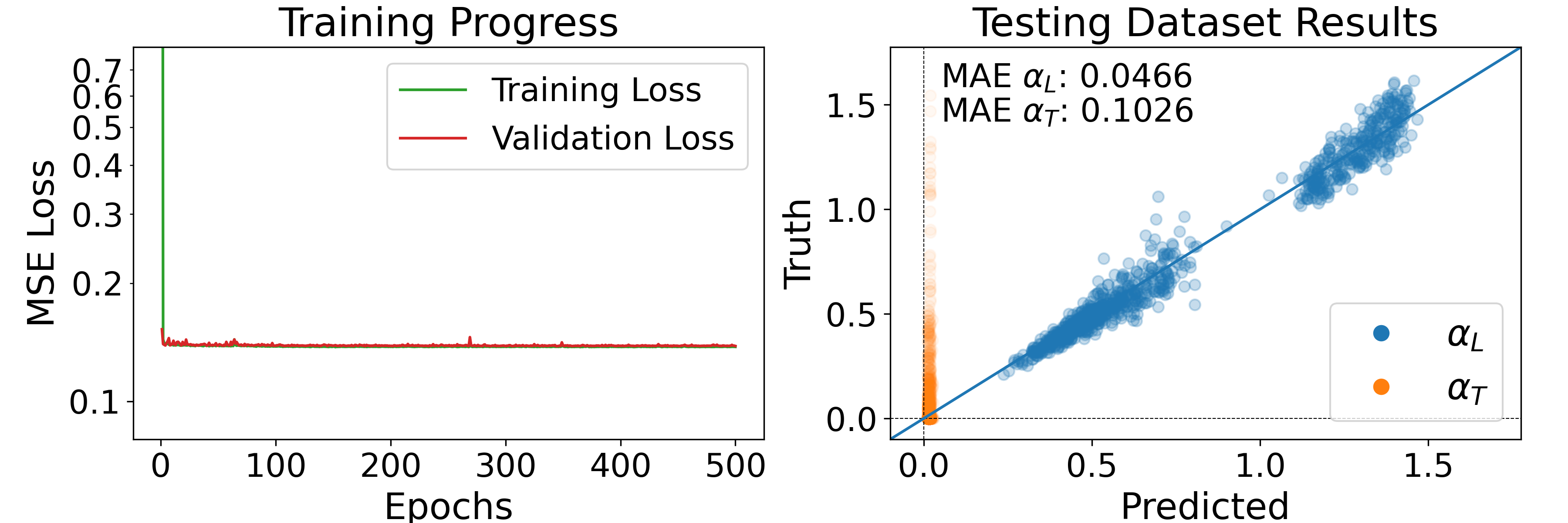}
	\caption{\label{fig:CNN_results} The performance of the dispersivity CNN is outlined. To the left, the evaluated loss, calculated for the training and a validation dataset at the end of each epoch, is shown over a series of epochs. This loss is based on the loss function, as outlined in \cref{eq:MSE}. To the right, the model's performance while evaluating on a testing dataset is shown, with $\boldsymbol{\alpha}$ evaluations shown in blue and orange. Here, the mean absolute error is shown (\cref{eq:MAE}).}
\end{figure}

During training, random batches of the training dataset are passed to the model, and the loss is evaluated. Working backwards, layer by layer, the gradient of the loss function in each layer is determined, and an update to each trainable parameter weight is added. This series of model predictions on a batched set, and back-propagation updates is performed on the whole training dataset to finish one epoch, and then is repeated for a series of 500 epochs. With increasing epochs, these losses decrease and converge to a final model state, as is seen in \cref{fig:CNN_results}. In this case, training converges after around 100 epochs, and the following converged epochs show that further training does not benefit the model. 

The testing dataset is also passed to the model in order to test the model with a fully independent dataset. In order to quantify the error with a metric easily comparable to the real data, the mean absolute error can be quantified. This is done with the following function:
\begin{equation}\label{eq:MAE}
\text{MAE} = \frac{1}{n} \sum_{i=1}^{n} \left| \boldsymbol{\alpha}_i - \hat{\boldsymbol{\alpha}}_i \right|, 
\end{equation}
where variables remain the same as those listed in \cref{eq:MSE}.

Here, it can be observed that a fairly accurate prediction is made for $\alpha_L$, with the data-driven model making predictions similar to those calculated from the pore-scale simulations. 

While predictions of $\alpha_T$ have an absolute error in a similar range as the error of $\alpha_L$, the true values are much smaller, making the presented error more significant. The data-driven model has difficulty predicting higher values for $\alpha_T$, hovering between $0.0$ and $0.05$. As discussed in \cref{subsec:OptimalAlphas}, the dispersion model chosen in this work, based on velocity and two length characteristics, is simple in comparison to other models with multiple parameters describing what the $\alpha_T$ parameter describes here. As these advanced dispersion models require further parameters, as well as a symmetry vector, this work has focused on the more commonly used dispersion model outlined here. With this in mind, the errors made when predicting $\alpha_T$ are more reasonable due to it's complex geometric definition and under-developed definition.
	
The image-based data-driven model offers a way to evaluate pore geometries without direct pore-scale analysis, but its complexity makes it difficult to interpret. With millions of weighted parameters in a deep non-linear structure, it is infeasible to understand how the model processes images to produce its regression output. To interpret the model, we rely on manual descriptions of the dataset. Many metrics exist to describe pore geometries, and analyzing their distribution and correlation with dispersivity values can help clarify the model’s function. These metrics are explored in \cref{sec:MetricsEvaluation,sec:MNN}.

\section{Metrics Based Evaluation}	\label{sec:MetricsEvaluation}

When working with data-driven models, it is quite important to develop an understanding of the parameter space upon which the model was trained. As these models often function poorly outside of this parameter space, one should try and maximize the training dataset's coverage of a possible application space. As this data-driven model works with pore-geometry images, descriptions of these pore-geometries must be taken to evaluate the training dataset's coverage. To address this, a series of physical metrics are introduced and collected from each pore geometry. These metrics can be categorized into three groups, void-solid interface metrics (\cref{subsec:fluidsolidmetrics}), volume ratio and distribution based metrics (\cref{subsec:volumeratio}), and connectivity based metrics (\cref{subsec:connectivitymetrics}). Here, the development of these metrics is outlined, their distribution in the dataset discussed, and their correlation to the dispersivity data investigated. 

\begin{figure}
	\centering
	\begin{subfigure}[t]{0.29\textwidth}
		\centering
		\includegraphics[width=0.785\textwidth]{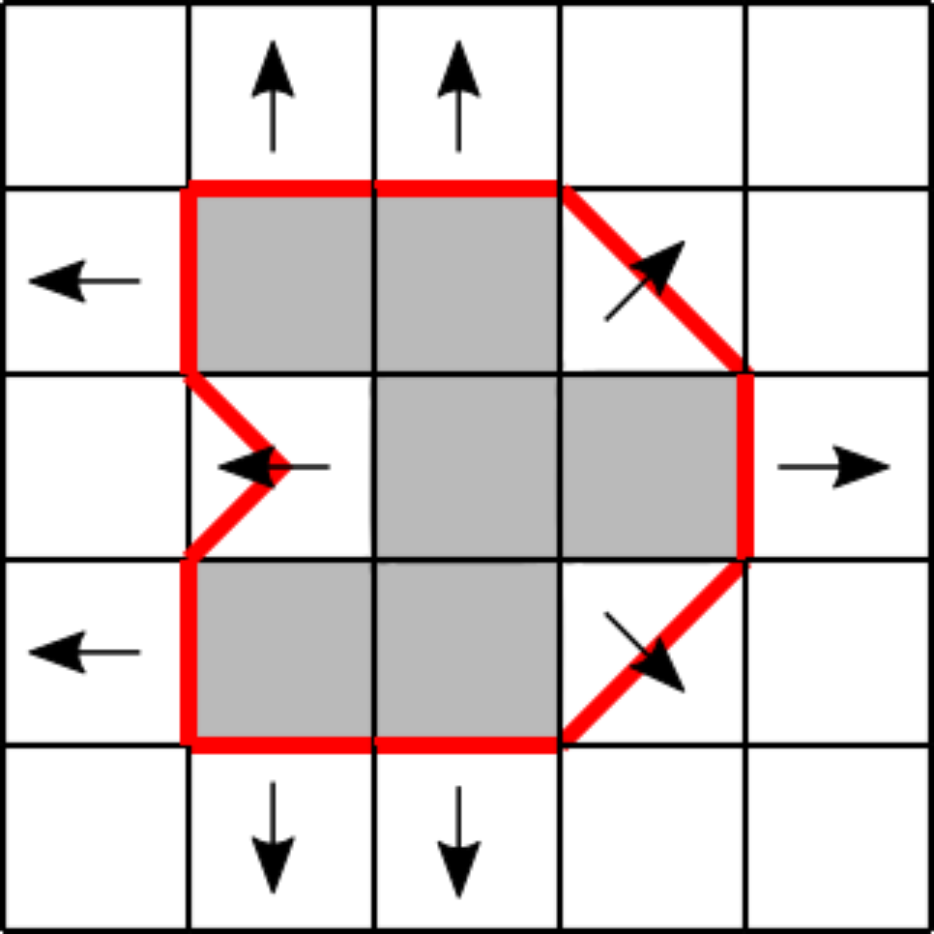}   
		\caption{An example determination of the surface properties. Pixels along void-solid interfaces, highlighted in red, contribute to the total surface area, with interface directionality indicated by arrows.}
		\label{fig:surface_example}
	\end{subfigure}
	\hspace{3mm}
	\begin{subfigure}[t]{0.29\textwidth}
		\centering
		\includegraphics[width=.80\textwidth]{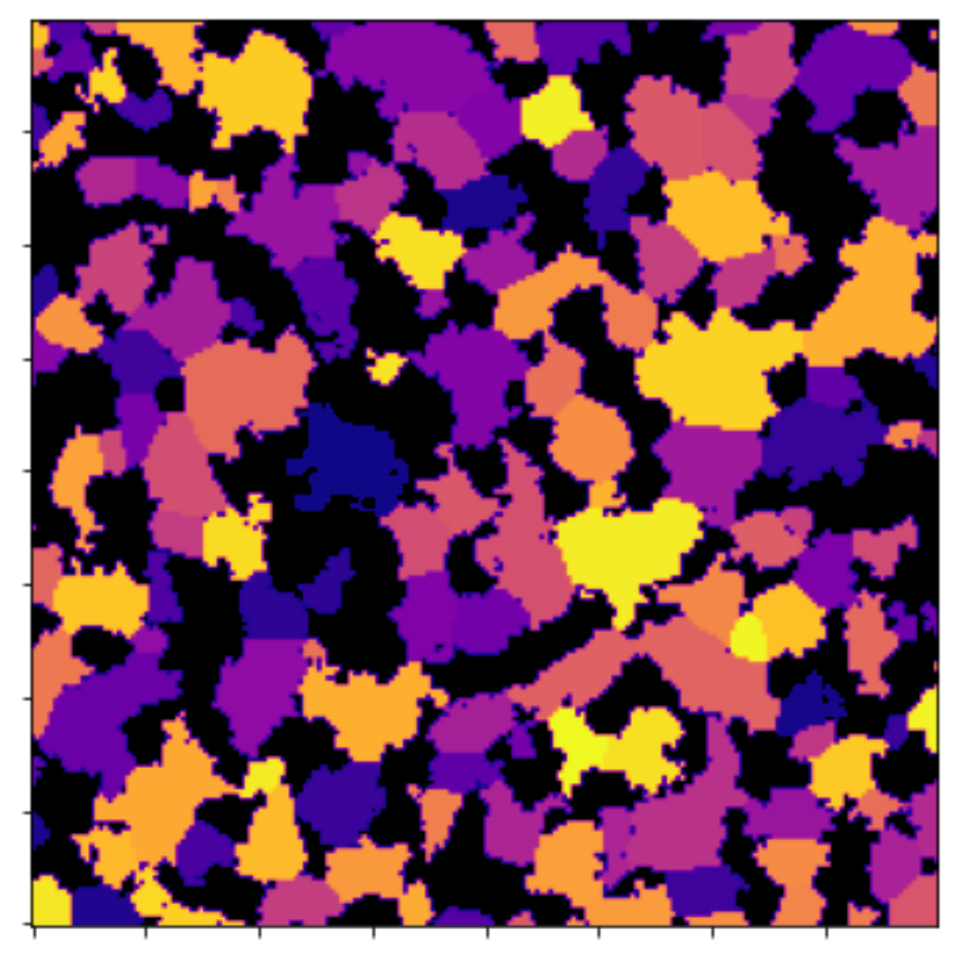}
		\caption{An example of pore geometry segmented into colored pore volumes. From these segments, parameters such as volume ratios and distribution can be calculated.}
		\label{fig:distribution_example}
	\end{subfigure}	
	\hspace{3mm}
	\begin{subfigure}[t]{0.29\textwidth}
		\centering
		\includegraphics[width=.80\textwidth]{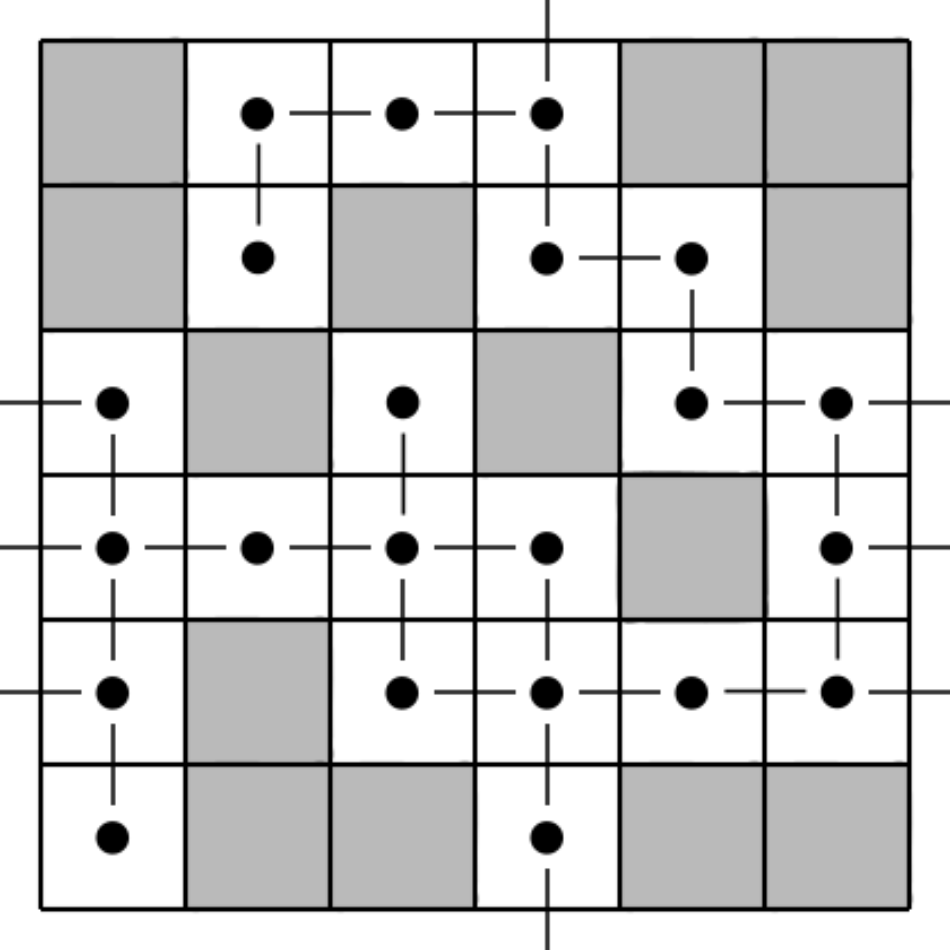}
		\caption{An example graph representation of a pore geometry. The graph nodes (black dots), fill the pore space pixels and are connected through edges (black connecting segments). }
		\label{fig:graph_scheme}
	\end{subfigure}
	\caption{\label{fig:metrics_graphics}}
\end{figure}

\subsection{Void-Solid Interface Metrics}\label{subsec:fluidsolidmetrics}

The pore geometry, consisting of domains $\Omega_{\text{solid}}$ and $\Omega_{\text{void}}$, can be characterized by the interface between them. Since flow conditions in porous media often depend on these void-solid interfaces, quantifying them helps describe pore geometries. This section outlines metrics for measuring the quantity, roughness, and directionality of these interfaces.

\subsubsection{Surface Area and Roughness}

The specific surface area, $\textbf{S}$ ($L^{-1}$), quantifies the total interface for each pore geometry, representing the ratio of the void-solid surface area to the total volume of the porous sample \cite{bear1988a}. The specific surface area is calculated using the algorithm shown in \cref{fig:surface_example}, where each pixel on the interface is evaluated. Depending on the configuration of pixel boundaries, different intersection lengths are added to the surface area. For only one pixel face on a surface, the intersection length is added to the surface area. If there are two pixel boundaries adjacent to each other, the diagonal pixel length is added, and if these two boundaries are opposite each other, both interfacial lengths are added, and in the case of three pixel boundaries, a concave surface is expected, and half of a diagonal length is added twice. 

For roughness, the fractal dimension, $R_{dim}$, can be calculated using a box-counting method and the following equation:

\begin{equation}
	R_{dim} = \lim_{L_{b} \rightarrow 0} \frac{-\log N_{b}(L_{b})}{\log (L_{b})}.
\end{equation}

Here, boxes of decreasing length, $L_b$, are placed over the image, stopping at a box sized at one pixel. The number of boxes containing both void and solid pixels at each length, $N_b(L_b)$, is then calculated. The resulting slope of their log relation then determines the roughness dimension, $R_{dim}$, where higher values correspond to rougher surfaces, and lower values to smoother surfaces \cite{gostick2019a}. 

\begin{figure}
	\centering
	\includegraphics[width=0.94\textwidth]{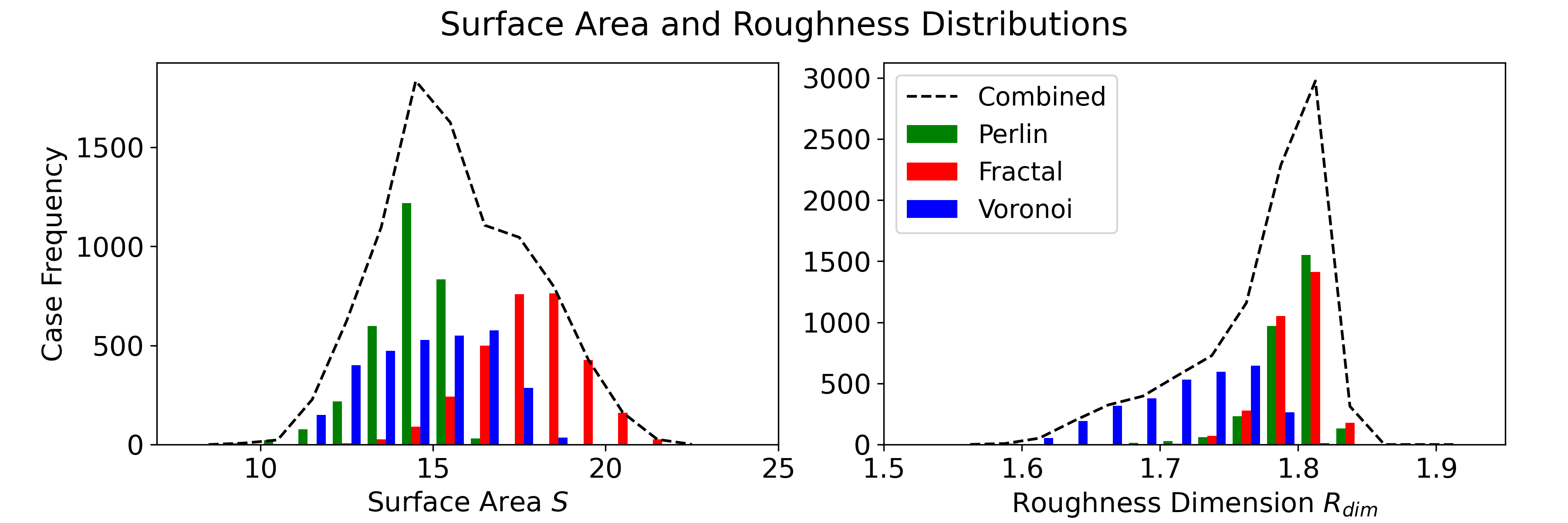}
	\caption{\label{fig:SurfaceAreaRoughnessDim} The distribution of surface areas,  $\textbf{S}$, and roughness dimensions, $R_{dim}$, found for the analyzed dataset.}
\end{figure}

As seen in \cref{fig:SurfaceAreaRoughnessDim}, the set of pore geometries used in the dataset vary in specific surface area, $\textbf{S}$, and roughness dimension, $R_{dim}$. The fractal cells generally have a higher surface area $\textbf{S}$, and the voronoi cells are generally smoother, with a lower roughness dimension. 

\subsubsection{Surface Directionality}

To describe the orientation of the void-solid interface, a directionality vector, $\mathbf{\Gamma}{d}$, is introduced, capturing each directional focus. Similar to calculating specific surface area, all interface pixels are analyzed, and the sum of their normal vectors is collected. The corresponding entry in $\mathbf{\Gamma}{d}$ is then incremented. \cref{fig:surface_example} shows examples of the normal surface vectors shown with arrows. 

Once all pixels are evaluated, the vector is normalized by scaling the surface area and pixel count. The normalized vector is analyzed using the standard deviation, $\sigma_{d}$ (deviation from average), skewness, $\gamma_{d}$ (asymmetry), and kurtosis, $\kappa_{d}$ (peaked bias from normal distribution).

\begin{figure}
	\centering
	\includegraphics[width=0.90\textwidth]{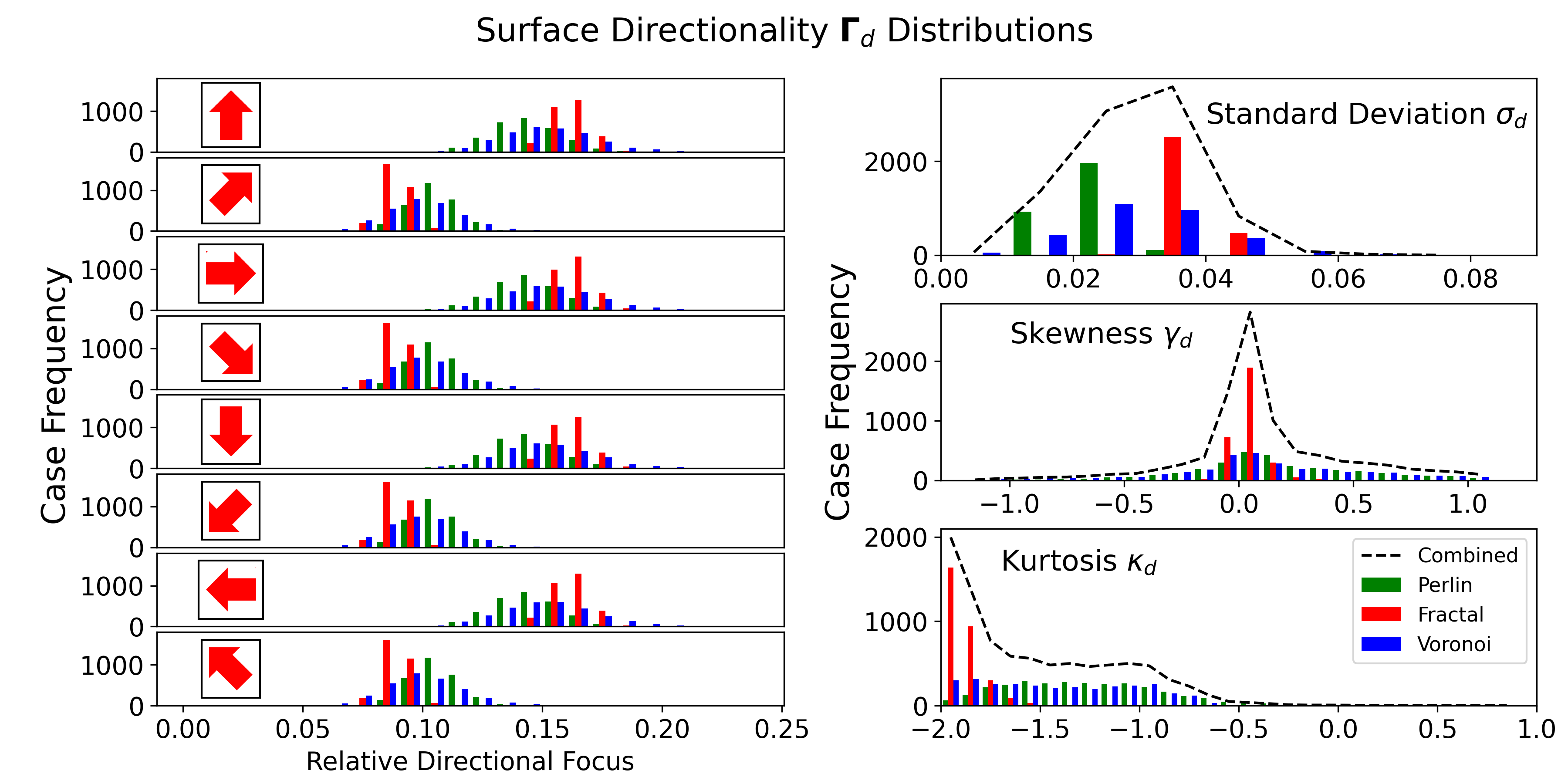}
	\caption{\label{fig:DirectionalityHists} The distributions of the relative directional focus in each direction, as well as the distribution of statistical descriptions of the directionality vector $\mathbf{\Gamma}_{d}$.}
\end{figure}

The relative directional focus varies for each direction, as shown in \cref{fig:DirectionalityHists}, where the theoretical average is 0.125. Without weighting the directional face length, primary axis directions occur more frequently, especially in fractal cases, which show a stronger preference due to their sharper, block-like roughness. Skewness from a symmetrical distribution is rare, with most values near 0. The kurtosis distribution confirms a bias toward axis-oriented faces.

\subsubsection{Dispersivity Correlation}

\begin{figure}
	\centering
	\includegraphics[width=0.96\textwidth]{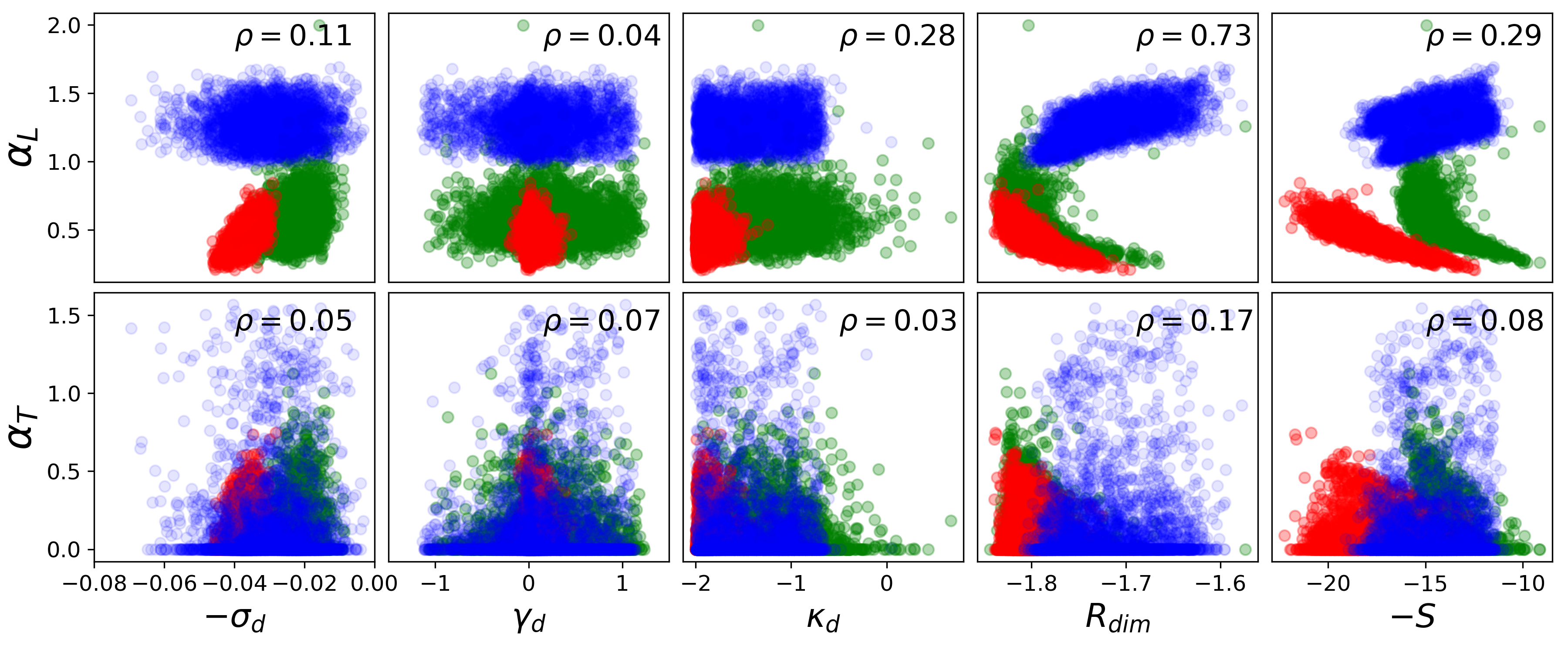}
	\caption{\label{fig:SurfaceCorrelation} The correlation between the surface-based metrics and the dispersivity parameters, $\alpha_L$ and $\alpha_T$. The $\rho$ shown in each box is the Pearson's coefficient of correlation.}
\end{figure}

As seen in \cref{fig:SurfaceCorrelation}, these surface descriptive metrics can then be compared against the dispersivity parameters, $\alpha_L$ and $\alpha_T$, in order to identify any possible correlations. 

Within each correlation plot, a total correlation value, Pearson's coefficient, $\rho$, is calculated according to the following equation:

\begin{equation}\label{eq:Pearsons}
\rho = \frac{\sum((x - \bar{x})(y - \bar{y}))}{\sqrt{\sum(x - \bar{x})^2 \sum(y - \bar{y})^2}},
\end{equation}

where $x$ is the evaluated metric and $y$ the dispersivity parameter and $\bar{x}$, and $\bar{y}$ their average values. This quantifies the strength of the linear correlation between the two, where values closer to 1 (or -1) indicate a higher correlation, and values closer to 0 point toward little linear correlation. 

As seen in \cref{fig:SurfaceCorrelation}, correlations between these surface-based metrics for $\alpha_L$ are generally higher than those for $\alpha_T$, although no strong linear correlation can directly be identified for either. Regarding $\alpha_L$, the roughness dimension ($R_{dim}$) correlates the highest, followed by the the surface area and the kurtosis of the directionality ($S$, $\kappa_{d}$). In addition, it should be noted, that when evaluated individually based on the geometry generation type, the correlations to $\alpha_L$ for roughness dimension and surface area increase significantly. All linear correlations for $\alpha_T$ score below $0.2$, meaning no clear linear relationship can be identified. 

\subsection{Volume Ratio and Distribution}\label{subsec:volumeratio}

One basic method for classifying porous media is by examining the volume ratios and relative pore sizes within a sample. Here, volume ratio and distribution-based metrics are developed, their dataset distribution is discussed, and their correlation with dispersivity parameters is evaluated.

\subsubsection{Pore Volume Ratios}

The volume ratio, or porosity $\phi$, is the ratio of void space $\Omega_{\text{void}}$ to the total volume $\Omega$: $\phi = \nicefrac{\Omega_{\text{void}}}{\Omega}$ \cite{bear1988a}. In this case, porosity is calculated as the total number of connected pore-space pixels divided by the total pixels in the image.

\subsubsection{Pore-Size Distribution}

\begin{figure}
	\centering
	\includegraphics[width=0.96\textwidth]{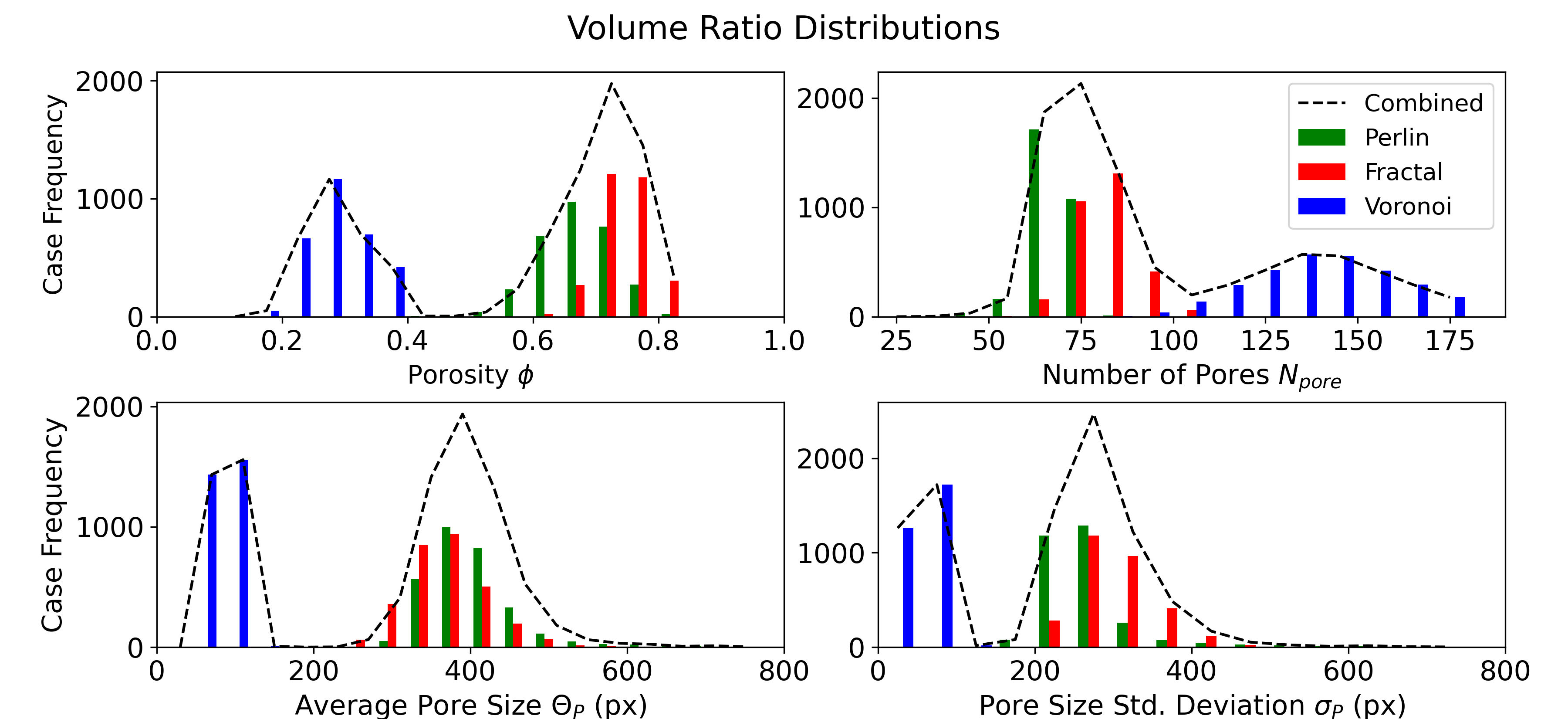}
	\caption{\label{fig:PoreSpaceHistograms} The porosity and pore distribution frequencies throughout the evaluated dataset.}
\end{figure}

In addition to volume ratio, pore-size distribution helps classify porosity in a pore geometry. Key metrics include the number of pores ($N_{pore}$), average pore size ($\Theta_{pore}$), and pore-size standard deviation ($\sigma_{pore}$). To calculate these, pore-network modeling techniques from \cite{gostick2017a, gostick2019a} are used, applying a segmentation algorithm to divide the pore space into individual pores, each with their own volume, as shown in \cref{fig:distribution_example}.

As seen in \cref{fig:PoreSpaceHistograms}, Perlin and fractal based geometries have higher porosities, larger average pore sizes, and greater variations when compared to Voronoi cells, which have smaller, more consistent pores. This is largely due to the higher porosity filters required to maintain connectivity across the pore geometry. Although these differences are clear, these ratios and parameters do not consider dead-end pore space, which is much more common among the higher porosity perlin and fractal cases. In addition, the segmentation of the regular pore spaces within voronoi cells can be difficult, and may be effected by the resolution of the image being processed. 
	
\subsubsection{Dispersivity Correlation}

\begin{figure}
	\centering
	\includegraphics[width=0.96\textwidth]{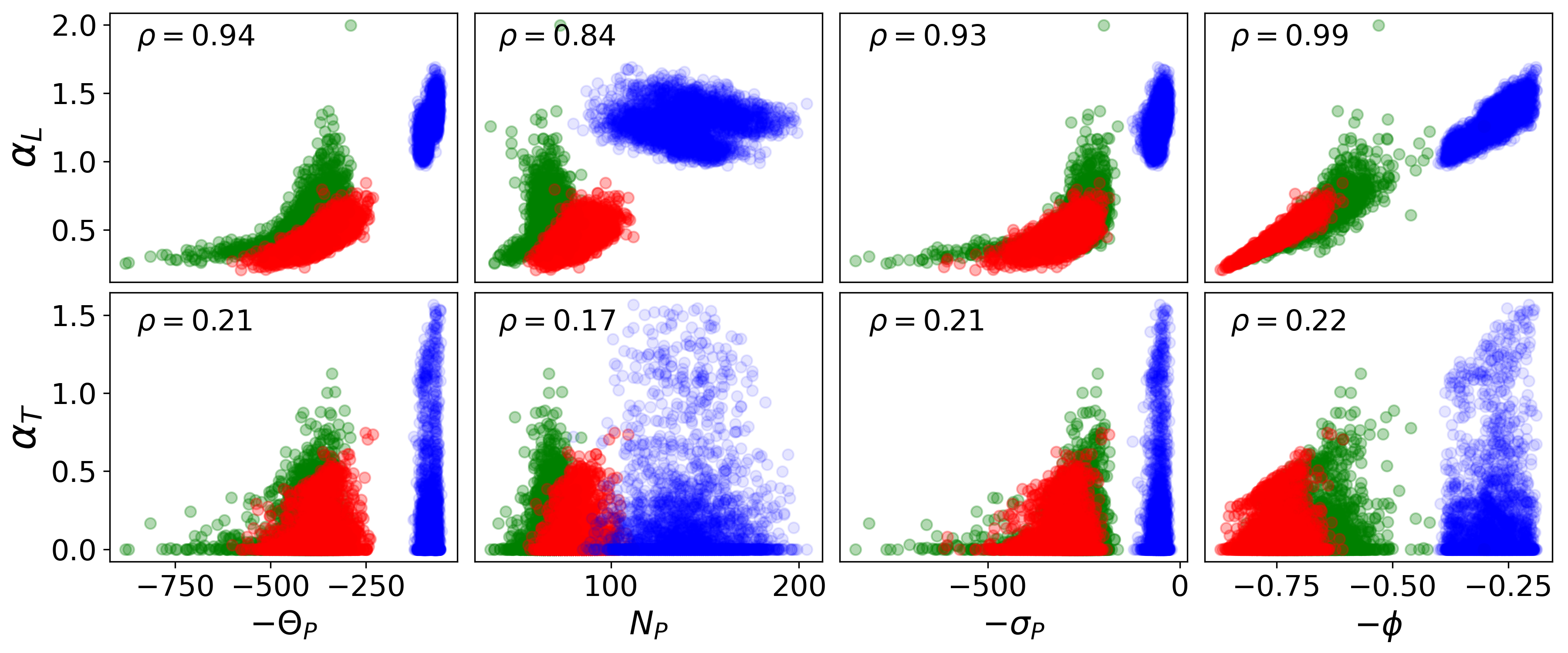}
	\caption{\label{fig:VolumeRatioCorrelation} The correlation between the volume ratio based metrics and the dispersivity parameters, $\alpha_L$ and $\alpha_T$. The $\rho$ shown in each box is the Pearson's coefficient of correlation.}
\end{figure}

As seen in \cref{fig:VolumeRatioCorrelation}, these volume ratio and distribution-based metrics can be compared with the dispersivity parameters, $\alpha_L$ and $\alpha_T$, and their linear correlation quantified with $\rho$, as defined in \cref{eq:Pearsons}. 

Again, correlations for $\alpha_L$ are all higher than the correlations outlined for $\alpha_T$. For $\alpha_L$, the porosity correlates strongly, followed by the average and standard deviation in pore size. All correlations for $\alpha_T$ are quite low, with porosity being the highest, and indicate no direct linear correlation. As noted previously, these correlations are developed for the full dataset. As the pore-space produced by each generator is quite different, isolated correlations my be higher. 

\subsection{Connectivity Focused Metrics}\label{subsec:connectivitymetrics}

In addition to metrics for the void-solid interface and pore-size distribution, pore connectivity can also be described. A common metric for this is permeability, either scalar or tensorial, which describes the geometric resistance to flow. Unfortunately, its exact determination from a pore-geometry is complex, requiring experiments or detailed simulations \cite{wagner2021}. Instead, we consider two alternative connectivity metrics important for flow and transport: tortuosity ($\tau$) and maximum flow ($F$).

Both metrics can be directly calculated from a pore geometry using a graphical representation of the pore space \cite{tang2021a, gaerttner2021a}, as shown in \cref{fig:graph_scheme}. Graph analysis, commonly used to evaluate network connectivity \cite{berge1985a}, is applied here by converting the pore-space domain $\Omega_{\text{void}}$ into a network graph. In this graph, each pixel in the pore space is represented by a node, and adjacent void space pixels are connected by edges. Traverses through the resulting network can then describe the shortest path between two nodes, or the general connectivity between various nodes. \cref{fig:graph_scheme} shows a pore geometry network graph where a pore geometry is represented as a network of connected nodes and edges. 

\begin{figure}
	\centering
	\begin{subfigure}{0.42\textwidth}
		\centering
		\includegraphics[width=0.9\textwidth]{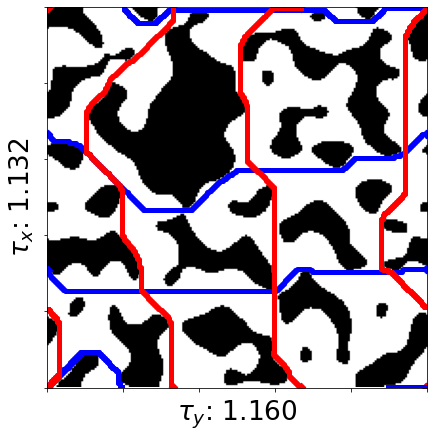}
		\subcaption{\label{fig:tortuosity_ex}}
	\end{subfigure}
	\begin{subfigure}{0.42\textwidth}
		\centering
		\includegraphics[width=0.9\textwidth]{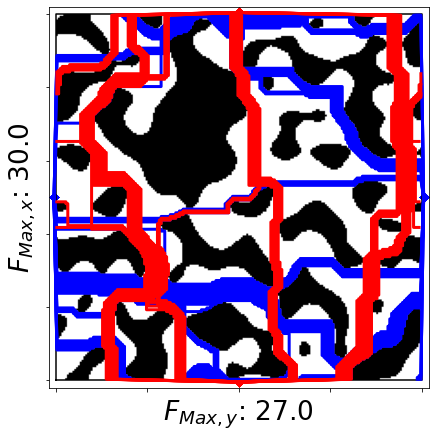}
		\subcaption{\label{fig:maxflow_ex}}
	\end{subfigure}
	\captionsetup{format=plain}
	\caption[Pore Metrics: Tortuosity and Maximum Flow]{\label{fig:graph_examples} Shown here are example tortuosity (\ref{fig:tortuosity_ex}) and maximum flow (\ref{fig:tortuosity_ex}) calculations. Paths in the x and y directions are shown in blue and red, respectively. These paths are used to trace the tortuous distances, as well as to find the edges to be removed for the maximum flow. Periodic connections are considered in each case. }
\end{figure}

\subsubsection{Tortuosity}

Tortuosity, $\tau$, was first introduced by \cite{carman1937a} to approximate the difference in flow path length between real chaotic porous media and the simplified bundle-of-tubes representation. Since then, it has been widely used in flow and transport models with varying definitions \cite{ghanbarian2013a}. In this work, geometrical tortuosity is used, defined as $\tau = \nicefrac{\left<T_{pore}\right>}{T_{L}}$, where $\left<T_{pore}\right>$ is the average geometric path length and $T_{L}$ is the straight-line distance.

To calculate $\left<T_{pore}\right>$, a geometric shortest path analysis is performed using a graph algorithm \cite{berge1985a, dijkstra1959a}. Sample points on opposite sides of the pore space are connected via the shortest path, and the average path length is divided by the straight-line distance to determine tortuosity. Directional tortuosities, $\tau_x$ and $\tau_y$, are also calculated by evaluating paths along each axis. This is illustrated on an example pore-geometry in \ref{fig:tortuosity_ex}. 

\begin{figure}
	\centering
	\includegraphics[width=0.96\textwidth]{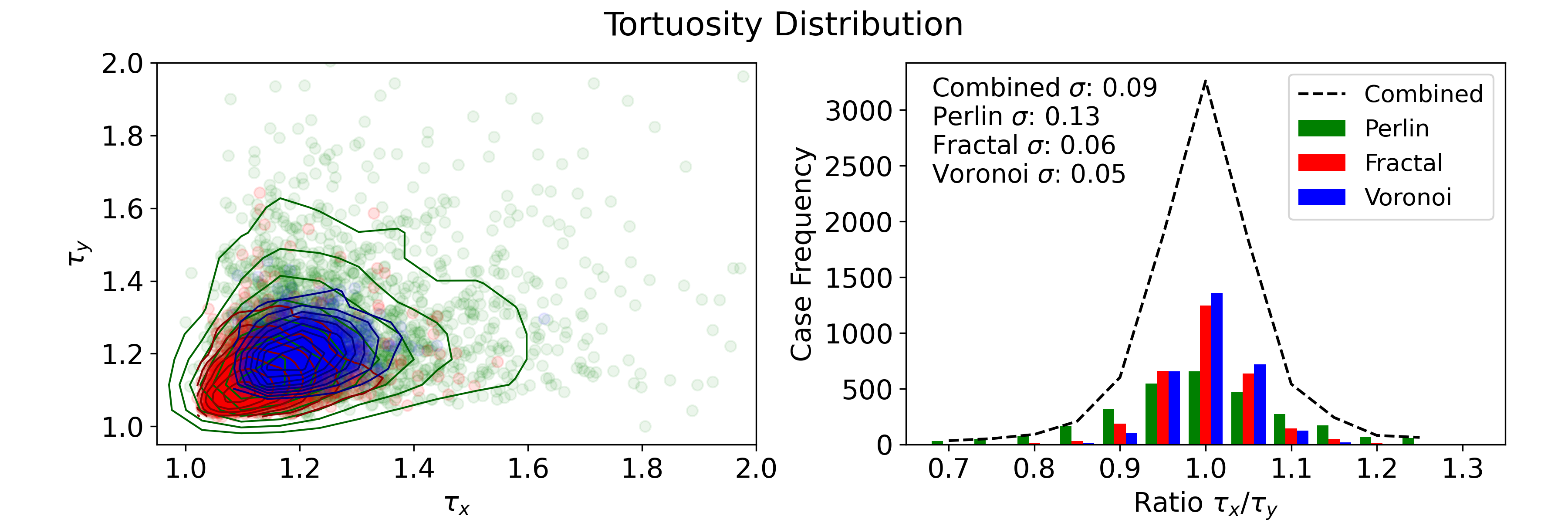}
	\caption{\label{fig:TortuosityHist} The distribution of the directional tortuosity $\tau_x$ and $\tau_y$ is shown to the left. The distribution of their ratio across the used dataset is shown to the right. To the right, the standard deviation of this distribution $\sigma$ is also shown.}
\end{figure}

As seen in \cref{fig:TortuosityHist}, the majority of the samples measure tortuosity between 1.0 and 1.4, in both the x and y directions. Further, the ratio of these directional tortuosities is evenly centered around a ratio of 1, with perlin cells producing a larger standard deviation. 

\subsubsection{Maximum Flow}

The maximum flow metric, $(F_x, F_y)$, reflects the size of bottlenecks in the porous medium and the geometric resistance to flow. Also known as the minimum cut, this metric represents the total number of pore-space edges that must be removed to disconnect the graph. This is done by iteratively removing connecting paths until none remain, with the cumulative removed edge weights determining the maximum flow \cite{fordfulkerson1956a}. A calculation of the maximum flow in each direction is shown for an example pore-geometry in \ref{fig:maxflow_ex}. 

Each pore geometry is analyzed to find the maximum flow in both directions $(F_x, F_y)$. Higher values indicate higher permeability \cite{gaerttner2021a, kanavas2021a}, while larger differences between $F_x$ and $F_y$ suggest flow anisotropy, with similar values indicating isotropy.

\begin{figure}
	\centering
	\includegraphics[width=0.96\textwidth]{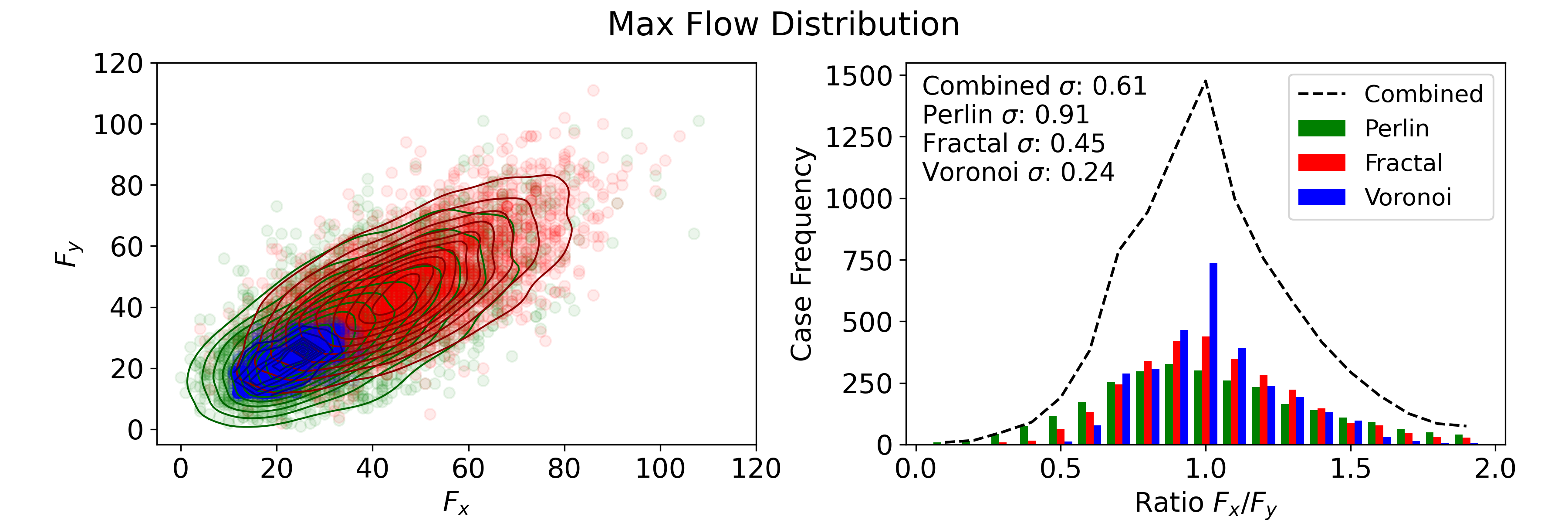}
	\caption{\label{fig:MaxmimumFlowHist} The distribution of the directional maximum flow $F_x$ and $F_y$ is shown to the left. The distribution of their ratio across the used dataset is shown to the right. To the right, the standard deviation of this distribution $\sigma$ is also shown.}
\end{figure}

As seen in \cref{fig:MaxmimumFlowHist}, voronoi cells cluster around values between 15 and 40, where perlin and fractal cells have a larger range, between one and 80. In addition, voronoi and fractal cells produce a tighter standard deviation in their ratio distribution, indicating fewer cells with larger anisotropy. 

\subsubsection{Dispersivity Correlation}

\begin{figure}
	\centering
	\includegraphics[width=0.96\textwidth]{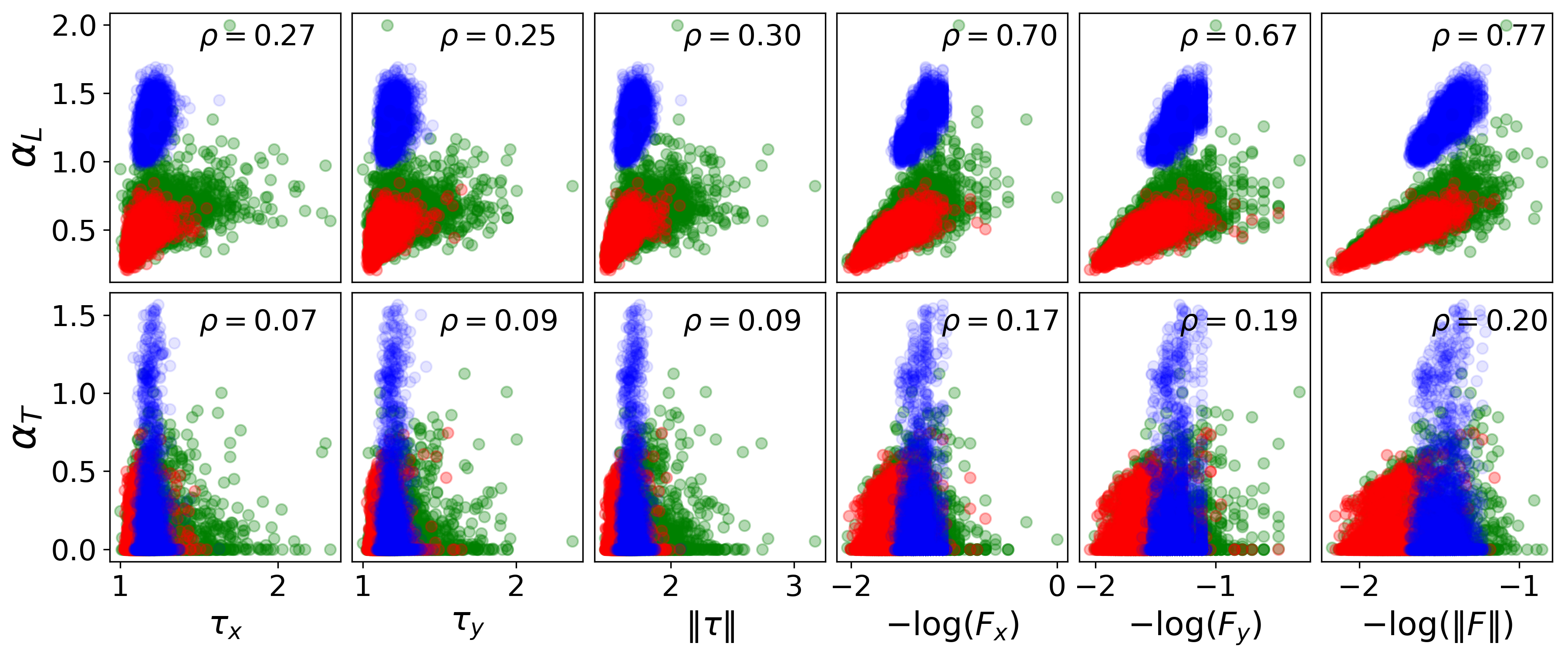}
	\caption{\label{fig:ConnectionsCorrelation} The correlation between the connectivity based metrics and the dispersivity parameters, $\alpha_L$ and $\alpha_T$. The $\rho$ shown in each box is the Pearson's coefficient of correlation.}
\end{figure}

As seen in \cref{fig:ConnectionsCorrelation}, these connectivity based metrics can then be compared against the dispersivity parameters, $\alpha_L$ and $\alpha_T$ and their linear correlation quantified with $\rho$, as defined in \cref{eq:Pearsons}. 

In this case, the logarithm of each directional maximum flow is shown individually, as well as the logarithm magnitude of each metric as a vector. Again, $\alpha_L$ consistently shows a stronger linear correlation with the connectivity metrics in comparison with the correlations for $\alpha_T$. The maximum flow parameter here shows the strongest linear correlations to the dispersivity parameters, $\alpha_L$ and $\alpha_T$. Again, correlations for $\alpha_T$ are much smaller when compared with $\alpha_L$. 

\section{Metrics Based Data-Driven Model}\label{sec:MNN}

To integrate the data-driven model from \cref{sec:datadrivenmodel} with the metrics-based analysis from \cref{sec:MetricsEvaluation}, one can replace the image input and processing steps in the CNN with a vector of descriptive metrics, followed by dense layers and an output. This approach offers a few benefits. If the performance matches that of the image-based model, it suggests the metrics sufficiently describe the pore geometry, covering the relevant parameter space. Additionally, this method simplifies the model significantly, replacing image-based inputs with scalar values and removing many unnecessary layers. 

Metrics-based neural networks, a form of multi-layer perceptron, has been used in other porous media studies to capture non-linear correlations between multiple variables to predict things like the permeability \cite{fu2023}, transport parameters \cite{taghizadeh2022a}, and multiphase parameters \cite{liu2019a}. While individual metrics may not show strong correlations with dispersivity parameters, their combination might be more effective. If this model performs similarly to the image-based one, it confirms that the metrics introduced are a good representation of the dataset.

\subsubsection{Model Architecture and Parameters} 

In order to maintain a structure comparable to the CNN outlined previously, the dense layers following the convolutional layers can be maintained, and only the image-processing layers removed. The metrics outlined in this case above are stacked to a vector of 21 data points (tortuosities, max flows, pore-size distribution statistics, porosities, roughness dimensions, surface areas, directionality vectors, and directionality statistics) and used as an input-data layer. As is done in the CNN, the dispersivity parameters, $\alpha_L$ and $\alpha_T$, are the target parameters set at the output layer of the model. Due to the reduced shape of the input layer, in comparison with all of the convolutional layers in the CNN, this model has much fewer parameters, with a total of only $p_{M} = ~34$ million. A summary of the layers used is shown in \cref{tab:MetricsNN}.

Before these data points are passed to the dense layers, they are normalized to a common scale to prevent any unwanted weighting. The loss associated with each dense layer is then regularized with a $\mathcal{L}^2$ regularization parameter of $1e-6$, in order to reduce any over-fitting as described in \cref{eq:MSE}. The same optimizer is used as in the CNN with a lower learning rate of $\epsilon_l = 2e-4$, and the same mean squared error loss function is used. In addition, the early stopping criterion is applied, such that the training stops when the validation loss ceases to decrease, again preventing over-fitting.
\begin{table}[h]
	\centering
	\begin{tabular}{|c|c|c|c|}
		\hline
		\textbf{Layer} & \textbf{Type} & \textbf{Shape} & \textbf{$\xi$ Count} \\		\hline
		1 & Input & 21 & 0 \\
		2 & Normalization & 21 & 43 \\ \hline
		3 & Dense Layer & 4096 & 90,112 \\
		4 & Dense Layer & 4096 & 16,781,312 \\
		$5_{a,b}$ & 2x Split Head Dense Layers & 2048 & 16,781,312 \\ \hline
		$6_{a,b}$ & 2x Regression Output & 1 & 4098 \\ \hline
	\end{tabular}
	\caption{\label{tab:MetricsNN} Metrics Based Neural Network}
\end{table}

\subsubsection{Metrics Based Model Performance} 

\begin{figure}
	\centering
	\includegraphics[width=0.99\textwidth]{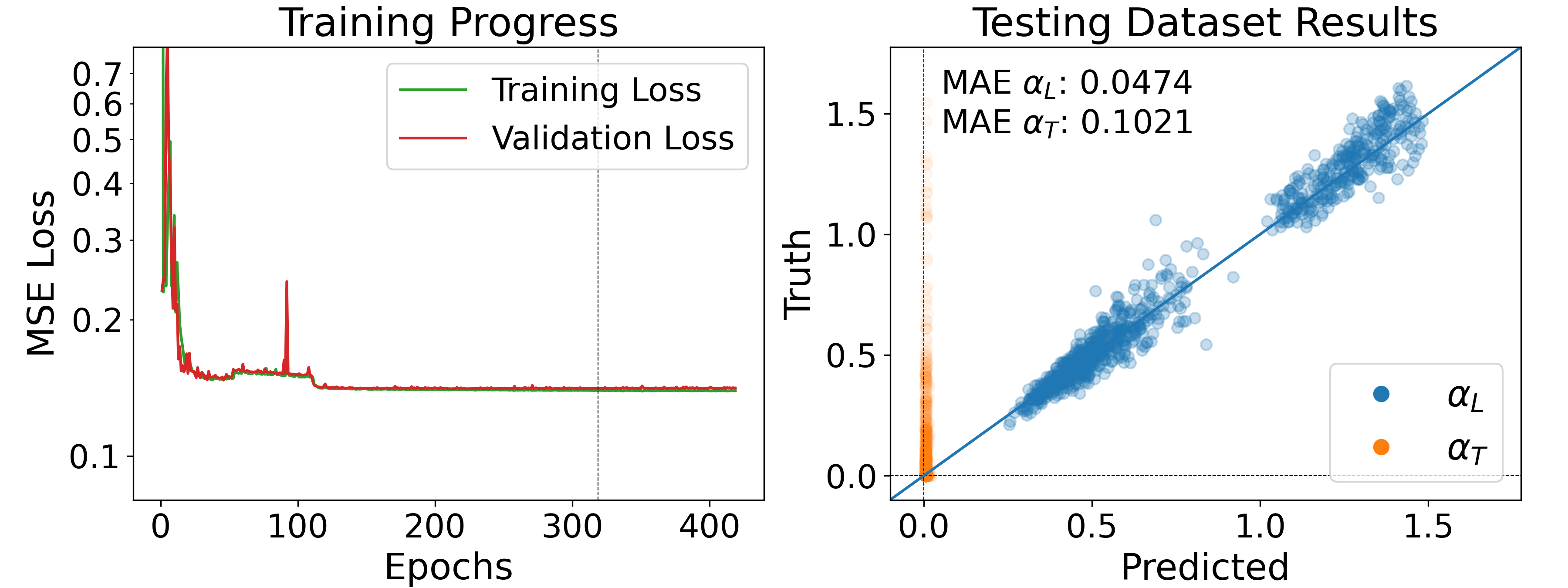}
	\caption{\label{fig:MetricsNN} The performance of the metrics based neural network is outlined. To the left, the evaluated total loss and individual losses as applied to the training and the validation dataset is shown over a series of epochs. Early stopping criteria are implemented such that the final model parameter state is collected where the dashed line is shown. To the right, the model's predictions on a testing dataset is shown. Each sample $\boldsymbol{\alpha}$ is shown, and the mean absolute error (\cref{eq:MAE}) for each variable is displayed.}
\end{figure}

As seen in \cref{fig:MetricsNN}, as the training epochs progress, the error as determined by the loss function, decreases to a point of convergence. According to the early stopping criteria, this convergence is identified and training is stopped to prevent any over-fitting. Training stops around epoch 415, as the validation error had not decreased for more than 100 epochs. 

When compared with the image-based CNN, the results are quite good for the longitudinal dispersivity parameter, with a mean absolute error of less than 0.5. The prediction of the transversal dispersivity is similarly poor when compared to the image-based model. With these models each producing similar results, one can assume that the metrics developed to describe the images are descriptive enough to train a model of reduced size to similar accuracy.  

As shown in the direct metrics-dispersivity correlations, $\alpha_L$ is more closely related to pore geometry. In contrast, predicting the transversal dispersivity, $\alpha_T$, has been more challenging, with larger errors in both data-driven models and weak correlations with descriptive metrics. Given the difficulties outlined in \cref{subsec:OptimalAlphas}, $\alpha_T$ likely under represents the true pore-scale transversal dynamics, making it difficult to directly extract for use in large-scale models based solely on pore geometry.

Another advantage of this metrics-based model over the image-based CNN is that it uses metrics that can be estimated for a porous medium. By estimating the distribution of these metrics, the model can predict a range of dispersivity values without needing detailed pore geometry images.

\section{Conclusions}	\label{sec:conclusions}
This work presents a data-driven framework for multi-scale parametrization of velocity-dependent dispersive transport in porous media. Pore-scale flow and transport simulations are performed on a series of periodic pore geometries. Next, the results are volume-averaged, isolating the dispersive transport, which is matched to optimal parameters at the REV scale through a dispersive closure term. These results are validated on simple unit cells with symmetric and asymmetric geometries for clarity. A convolutional neural network is then trained to predict these dispersivity parameters based on images of the pore geometries, replacing the pore-scale simulations with a data-driven model. Finally, descriptive metrics are introduced to better understand the parameter space and establish correlations between dispersivity and other REV-scale metrics. These metrics are then used to build a neural network that predicts dispersivity based on the descriptive metrics alone.

The longitudinal dispersivity parameter, $\alpha_L$, is accurately predicted from descriptions of the pore geometry, with good results from both image-based and metrics-based neural networks, and stronger correlations to descriptive metrics. In contrast, the transverse dispersivity, $\alpha_T$, is poorly predicted and shows limited correlations to all descriptive metrics. This is likely because the dispersion model used here is not always able to describe pore-scale transversal dispersivity. While $\alpha_L$ can be predicted from a pore geometry, $\alpha_T$ may require a more complex model that accounts for different geometries and symmetries \cite{lichtner2002a}.

Extensions of this analysis can be explored in several ways. First, alternative dispersivity models, from simpler single-parameter versions to more complex ones, could be investigated. Additional parameters would need to be balanced against model accuracy, and any correlation with the pore metrics introduced here. A similar approach could also be used to parametrize drag terms for porous flow with inertial effects or assess how flow and transport parameters change with pore-space variations. Finally, applying symbolic regression to the descriptive metrics correlations could further improve the interpretability of the analysis, and help to develop more useful heuristic models.

\vspace{-2mm}
\section*{Acknowledgments}
\vspace{-2mm}
The authors would like to thank Prof. Dr.-Ing Andrea Beck for her guidance and constructive discussions. 
All authors would like to thank the German Research Foundation (DFG) for supporting this work by funding SimTech via Germany’s Excellence Strategy (EXC 2075 – 390740016).
The authors would also like to thank the DFG for supporting this work by funding SFB 1313, Project Number 327154368, Research Project A02.
\vspace{-5mm}
\section*{Declarations}
\vspace{-2mm}
The authors do not have any relevant financial or non-financial affiliations to declare.

\clearpage
\bibliography{literature}
	
\end{document}